\documentclass[reqno,a4paper]{amsart}

\usepackage{amsmath}
\usepackage{amsthm}
\usepackage{amssymb}
\usepackage{mathrsfs}
\usepackage[latin1]{inputenc}
\usepackage{diagrams}
\usepackage{a4wide}

\DeclareMathOperator{\re}{re}
\DeclareMathOperator{\im}{im}

\newcommand{\half}{{\textstyle\frac{1}{2}}}
\newcommand{\third}{{\textstyle\frac{1}{3}}}
\newcommand{\fourthird}{{\textstyle\frac{4}{3}}}
\newcommand{\quart}{{\textstyle\frac{1}{4}}}
\newcommand{\sixth}{{\textstyle\frac{1}{6}}}
\newcommand{\bbz}{\mathbb{Z}}
\newcommand{\bbr}{\mathbb{R}}
\newcommand{\bbc}{\mathbb{C}}
\newcommand{\bbo}{\mathbb{O}}

\newcommand{\bbrp}{\mathbb{R}^{+}}
\newcommand{\gstr}{$G_{2}$-structure}
\newcommand{\cystr}{Calabi-Yau structure}
\newcommand{\slstr}{$SL(\bbc^{3})$-structure}
\newcommand{\gmfd}{$G_{2}$-manifold}

\newcommand{\ltwoorth}{$L^{2}$-orthogonal}

\newcommand{\sob}[1]{L^{2}_{#1}}

\newcommand{\sobx}[2]{L^{2}_{#1,#2}}
\newcommand{\holda}[1]{C^{#1,\alpha}}
\newcommand{\holdax}[2]{C^{#1,\alpha}_{#2}}
\newcommand{\holdad}[1]{C^{#1,\alpha}_{\delta}}
\newcommand{\holdadasy}[1]{$\holdad{#1}$\nobreakdash-\hspace{0pt}asymptotic}

\newcommand{\harm}{\mathcal{H}}
\newcommand{\norm}[1]{\Vert #1 \Vert}

\newcommand{\lgnorm}[1]{\Vert #1 \Vert_{L^{2}(g)}}
\newcommand{\lgpnorm}[1]{\Vert #1 \Vert_{L^{2}(g')}}

\newcommand{\calm}[2]{\mathcal{#1}^{#2}}

\newcommand{\calq}[1]{\mathcal{#1}_{\mathcal{Q}}}
\newcommand{\calp}[1]{\mathcal{#1}_{\scriptscriptstyle+}}
\newcommand{\calx}[1]{\mathcal{#1}_{X}}
\newcommand{\calxa}{\mathcal{X}'_{\delta}}
\newcommand{\cald}[1]{\mathcal{#1}_{\delta}}

\newcommand{\orbt}{T}
\newcommand{\slt}{K}
\newcommand{\dkaq}{\mathcal{D}^{k+1,\alpha}_{\mathcal{Q}}}

\newcommand{\ztwo}{\mathcal{Z}^{2}}
\newcommand{\ztre}{\mathcal{Z}^{3}}
\newcommand{\ztreq}{\mathcal{Z}^{3}_{\mathcal{Q}}}
\newcommand{\zfive}{\mathcal{Z}^{5}}
\newcommand{\inner}[1]{<\!\!#1\!\!>}
\newcommand{\linner}[1]{<\!\!#1\!\!>_{L^{2}}}
\newcommand{\lginner}[1]{<\!\!#1\!\!>_{L^{2}(g)}}
\newcommand{\lgpinner}[1]{<\!\!#1\!\!>_{L^{2}(g')}}

\newcommand{\contrat}{\textstyle\frac{\partial}{\partial t}}

\newcommand{\tildm}{\tilde M}
\newcommand{\inc}{e}
\newcommand{\hcup}{\cup}
\newcommand{\dirac}{\eth}

\numberwithin{equation}{section}
\newtheorem{thm}{Theorem}[section]
\newtheorem{prop}[thm]{Proposition}
\newtheorem{lem}[thm]{Lemma}

\newtheorem{cor}[thm]{Corollary}
\theoremstyle{definition}
\newtheorem{defn}[thm]{Definition}
\theoremstyle{remark}
\newtheorem{rmk}[thm]{Remark}
\newtheorem{ex}[thm]{Example}

\begin{document}

\title{Deformations of asymptotically cylindrical $G_{2}$-manifolds}
\author{Johannes Nordström}

\address{DPMMS, University of Cambridge,
    \mbox{Wilberforce} Road, Cambridge CB3 0WB, UK}
\email{j.nordstrom@dpmms.cam.ac.uk}

\begin{abstract}
We prove that for a $7$-dimensional manifold $M$ with cylindrical ends
the moduli space of exponentially asymptotically cylindrical
torsion-free \gstr s is a smooth manifold (if non-empty),
and study some of its local properties.
We also show that the holonomy of the induced metric
of an exponentially asymptotically cylindrical \gmfd{} is exactly $G_{2}$
if and only if the fundamental group $\pi_{1}(M)$ is finite
and neither $M$ nor any double cover of $M$ is homeomorphic to a cylinder.
\end{abstract}

\maketitle

\section{Introduction}
The holonomy group $G_{2}$ appears as an exceptional case in Berger's
classification of the Riemannian holonomy groups \cite{berger55}.
A metric with holonomy contained in $G_{2}$ can be defined
in terms of a parallel non-degenerate differential $3$-form,
often called a torsion-free \gstr.
A manifold $M$ has cylindrical ends
if a complement of a compact subset of $M$ is identified with
$X \times \bbrp$, for some compact manifold $X$ called the cross-section
of $M$. An important class of metrics on such manifolds are the
exponentially asymptotically cylindrical (EAC) ones, i.e. metrics which are
exponentially asymptotic to a product metric on the cylindrical part.
An EAC \gmfd{} is a $7$-dimensional EAC manifold $M$ with holonomy contained
in $G_{2}$, and its metric is defined by a torsion-free EAC \gstr.
Section $2$ covers this background in more detail.

On an EAC \gmfd{} $M^{7}$ the group of EAC diffeomorphisms isotopic to the
identity acts on the space
of torsion-free EAC \gstr s by pull-backs. We define the moduli space
of torsion-free EAC \gstr s to be the resulting quotient $\calp{M}$.
The precise definition involves a normalisation, which can also be interpreted
as dividing by the rescaling action of $\bbrp$ (see remark \ref{normalisermk}).
The main result of the paper is theorem \ref{maincylthm}, which states
that on an EAC \gmfd{} $M^{7}$
the moduli space $\calp{M}$ is a smooth manifold.
The proof of theorem \ref{maincylthm} is a generalisation of an
argument for the compact case outlined by Hitchin in~\cite{hitchin00}.
If $X^{6}$ is the cross-section of $M$ then
the dimension of the moduli space is given by the formula
\[ \dim \calp{M} = b^{4}(M) + \half b^{3}(X) - b^{1}(M) - 1 . \]

The cross-section $X$ of a manifold $M$ with cylindrical ends can be regarded
as the `boundary at infinity' of $M$, in the sense that $M$ can be identified
with the interior of a compact manifold with boundary $X$.
The asymptotic limit of an
EAC torsion-free \gstr{} on $M$ induces a Calabi-Yau structure on $X$, and
the proof of theorem \ref{maincylthm} requires understanding of the
deformations of this structure on the boundary. Theorem \ref{cymodthm}
states that the moduli space $\mathcal{N}$ of Calabi-Yau structures
on a compact connected oriented manifold $X^{6}$ is a smooth manifold of
dimension $b^{3}(X) + b^{2}(X) - b^{1}(X) - 1$. This is a special case of a more
general result due to Tian \cite{tian87} and Todorov \cite{todorov89}.
The argument given here is based on an elementary application of the
implicit function theorem, and is helpful for proving theorem \ref{maincylthm}.

There is natural `boundary map' $B : \calp{M} \to \mathcal{N}$ which sends
a class of torsion-free \gstr s on $M$ to the class of \cystr s that their
asymptotic limit defines on $X$.
Theorem \ref{modsubmersethm} states that this map is
a submersion onto its image, which is a submanifold of $\mathcal{N}$, and
an open subset of a subspace $\mathcal{N}_{A}$
determined by the topology of the pair $(M,X)$.

The final result of the paper is a
topological criterion for when the holonomy of the metric associated to
a torsion-free EAC \gstr{} is exactly $G_{2}$, rather than a subgroup.
Theorem \ref{maintopcondthm} states that $Hol(M) = G_{2}$
for an EAC \gmfd{} $M$
if and only if the fundamental group $\pi_{1}(M)$ is finite and
neither $M$ nor any double cover of $M$ is homeomorphic to a cylinder.

Precise statements for all the main results are given in section \ref{statesec}.
Section \ref{cymodsec} provides the needed deformation theory for compact
Calabi-Yau $3$-folds, extending a construction of the moduli space of
torsion-free \slstr s by Hitchin in \cite{hitchin00}.
In section \ref{hodgesec} we review Hodge theory for EAC manifolds,
which is an important tool in the proofs of the main results.
The main theorem \ref{maincylthm}
is proved in section $6$ along with results about the local properties
of the moduli space, and theorem \ref{maintopcondthm} is proved in section $7$.

It is complicated to produce examples of \gmfd s with holonomy exactly $G_{2}$.
The first compact examples were constructed by Joyce in~\cite{joyce96-I}.
Also, in \cite{kovalev03} Kovalev constructs non-trivial EAC \gmfd s
(with holonomy $SU(3)$), and uses them in a gluing construction to
produce compact manifolds with holonomy exactly $G_{2}$.
A future paper \cite{jn2} will show how to construct EAC manifolds with
holonomy exactly $G_{2}$, adapting the methods of \cite{joyce96-I} to the EAC
context.

\section{Background material}

In this section we review definitions and notation for Riemannian holonomy,
\gstr s on $7$-manifolds, \cystr s on $6$-manifolds, and manifolds with
cylindrical ends.

\subsection{Holonomy}
\label{holsub}

We define the holonomy group of a Riemannian manifold.
For a fuller discussion of holonomy see e.g. Joyce \cite[Chapter $2$]{joyce00}.

\begin{defn}
Let $M^{n}$ be a manifold with a Riemannian metric $g$.
If $x \in M$ and $\gamma$ is a closed piecewise $C^{1}$ loop in $M$ based at
$x$ then the parallel transport around~$\gamma$
(with respect to the Levi-Civita connection of the metric)
defines an orthogonal linear map
$P_{\gamma} : T_{x}M \to T_{x}M$. The \emph{holonomy group}
$Hol(g, x) \subseteq O(T_{x}M)$ at $x$ is the group generated
by $\{P_{\gamma} : \gamma \textrm{ is a closed loop based at } x\}$.
\end{defn}

If $x, y \in M$ and $\tau$ is a path from $x$ to $y$ we can define a
group isomorphism $Hol(g, x) \to Hol(g, y)$ by 
$P_{\gamma} \mapsto P_{\tau} \circ P_{\gamma} \circ P^{-1}_{\tau}$.
Provided that $M$ is connected we can therefore identify $Hol(g, x)$ with
a subgroup of $O(n)$, independently of $x$ up to conjugacy,
and talk simply of the holonomy group of~$g$.

There is a correspondence between tensors fixed by the holonomy group and
parallel tensor fields on the manifold.

\begin{prop}[{\cite[Proposition $2.5.2$]{joyce00}}]
\label{simpleholprop}
Let $M^{n}$ be a manifold with Riemannian metric $g$, $x \in M$
and $E$ a vector bundle on $M$ associated to $TM$.
If $s$ is a parallel section of $E$ then $s(x)$ is fixed by
$Hol(g, x)$.
Conversely if $s_{0} \in E_{x}$ is fixed by $Hol(g, x)$ then there is a
parallel section $s$ of $E$ such that $s(x) = s_{0}$.
\end{prop}

\subsection{$G_{2}$-structures}

An effective approach to \gstr s is to define them in terms
of \emph{stable $3$-forms}. Here we outline the properties of \gstr s, and
explain their relation to metrics with holonomy $G_{2}$. For a more complete
explanation see e.g. \cite[Chapter $10$]{joyce00}.

Recall that $G_{2}$ can be defined as the automorphism group of the normed
algebra of octonions. Equivalently, $G_{2}$ is the stabiliser in
$GL(\bbr^{7})$ of
\begin{equation}
\label{g2formeq}
\varphi_{0} = dx^{123} + dx^{145} + dx^{167} + dx^{246}
 - dx^{257} - dx^{347} - dx^{356} \in \Lambda^{3}(\bbr^{7})^{*} .
\end{equation}
If $V$ is a dimension $7$ real vector space and $\varphi \in \Lambda^{3}V^{*}$
then we call $\varphi$ \emph{stable} if it is equivalent to $\varphi_{0}$ under
some isomorphism $V \cong \bbr^{7}$.
We denote the set of stable 3-forms by $\Lambda^{3}_{+}V^{*}$.
Since the action of $GL(V)$ on $\Lambda^{3}V^{*}$
has stabiliser $G_{2}$ at a stable form $\varphi$, and $\dim G_{2} = 14$,
it follows by dimension-counting that $\Lambda^{3}_{+}V^{*}$ is open in
$\Lambda^{3}V^{*}$.
Since $G_{2} \subset SO(7)$ each $\varphi \in \Lambda^{3}_{+}V^{*}$ naturally
defines an inner product $g_{\varphi}$ and an orientation.

\begin{defn}
Let $M^{7}$ be an oriented manifold.
A \emph{\gstr} on $M$ is a section $\varphi$
of $\Lambda^{3}_{+}T^{*}M$ which defines the given orientation on $M$.
\end{defn}

Since a \gstr{} $\varphi$ on $M$ induces a Riemannian metric $g_{\varphi}$
it also defines a Levi-Civita connection $\nabla_{\varphi}$,
a Hodge star $*_{\varphi}$ and a codifferential~$d^{*}_{\varphi}$.

\begin{defn}
\label{torsionfreedef}
A \gstr{} $\varphi$ on an oriented manifold $M^{7}$ is \emph{torsion-free}
if \mbox{$\nabla_{\varphi} \varphi = 0$}.
A~\emph{\gmfd} is an oriented manifold $M^{7}$ equipped with a torsion-free
\gstr{} $\varphi$ and the associated Riemannian metric $g_{\varphi}$.
\end{defn}

Gray observed that
a \gstr{} is torsion-free if and only if it is closed and coclosed.

\begin{thm}[{\cite[Lemma $11.5$]{salamon89}}]
\label{graythm}
A \gstr{} $\varphi$ on $M^{7}$ is torsion-free if and only if
$d\varphi = 0$ and $d^{*}_{\varphi}\varphi = 0$.
\end{thm}

As an immediate application of proposition \ref{simpleholprop}
we have that metrics with holonomy contained in $G_{2}$ correspond
to torsion-free \gstr s.

\begin{cor}
\label{g2holcor}
Let $M^{7}$ be a manifold with Riemannian metric $g$. Then $Hol(g)$ is a
subgroup of $G_{2} \subset O(7)$ if and only if there is a torsion-free \gstr{}
$\varphi$ on $M$ such that $g = g_{\varphi}$.
\end{cor}

The condition that $Hol(g) \subseteq G_{2}$ imposes algebraic constraints on
the curvature of $g$. In particular

\begin{thm}[{\cite[Proposition $11.8$]{salamon89}}]
The metric of a torsion-free \gstr{} $\varphi$ is Ricci-flat.
\end{thm}

The deformation problem for torsion-free \gstr s on a compact oriented manifold
$M^{7}$ was solved by Joyce in \cite{joyce96-I}. Let $\mathcal{X}$ be
the space of smooth torsion-free \gstr s on $M$, and $\mathcal{D}$ the group of
diffeomorphisms of $M$ isotopic to the identity.
$\mathcal{D}$ acts on $\mathcal{X}$ by pull-backs, and
the \emph{moduli space of torsion-free \gstr s} on $M$ is the space of orbits
$\mathcal{M} = \mathcal{X}/\mathcal{D}$. Elements of $\mathcal{X}$
are closed $3$-forms, so define cohomology classes. This gives a well-defined
map to de Rham cohomology $\pi_{H}: \mathcal{M} \to H^{3}(M)$.
Joyce proved that

\begin{thm}[{\cite[Theorem $10.4.4$]{joyce00}}]
\label{maincptthm}
Let $M^{7}$ be a compact \gmfd. Then $\mathcal{M}$ is a smooth manifold
of dimension $b^{3}(M)$, and $\pi_{H} : \mathcal{M} \to H^{3}(M)$
is a local diffeomorphism.
\end{thm}

The main result of this paper generalises theorem \ref{maincptthm}
to the case when $M$ is an EAC \gmfd.

\subsection{Calabi-Yau $3$-folds}
\label{cydefsub}

One common definition of a Calabi-Yau manifold is that it is a Riemannian
manifold $X^{2n}$ with holonomy contained in $SU(n)$.
The stabiliser in $G_{2}$ of a vector in $\bbr^{7}$ is isomorphic to $SU(3)$,
and Calabi-Yau $3$-folds will appear naturally as the cross-sections of
EAC \gmfd s.

For our purposes it is convenient to define a \emph{\cystr} on a
$6$-dimensional manifold $X$ in terms of a pair of closed differential forms
$(\Omega,\omega)$.
This will make the relation to \gstr s clear. The role of the `stable'
$3$-form $\Omega$ is discussed by Hitchin in \cite{hitchin00}.
Let 
\begin{subequations}
\label{su3formseq}
\begin{gather}
\label{sl3formeq}
\Omega_{0} = dx^{135} - dx^{146} - dx^{236} - dx^{245} \in
\Lambda^{3}(\bbr^{6})^{*}, \\
\label{kahlerformeq}
\omega_{0} = dx^{12} + dx^{34} + dx^{56} \in \Lambda^{2}(\bbr^{6})^{*} .
\end{gather}
\end{subequations}
For an oriented real vector space $V$ of dimension $6$
let $\Lambda^{3}_{+}V^{*}$ be the set of $\Omega \in \Lambda^{3}V^{*}$
such that $\Omega$ is equivalent to $\Omega_{0}$
under some linear isomorphism $V \cong \bbr^{6}$.
We call such $\Omega$ \emph{stable}.
If we identify $\bbr^{6}$ with $\bbc^{3}$ by taking $z^{1} = x^{1} + ix^{2}$,
$z^{2} = x^{3} + ix^{4}$, $z^{3} = x^{5} + ix^{6}$ then
\begin{subequations}
\begin{gather*}
\Omega_{0} = \re \; dz^{1} \wedge dz^{2} \wedge dz^{3} , \\
\omega_{0} = {\textstyle\frac{i}{2}}(dz^{1} \wedge d \bar z^{1} +
dz^{2} \wedge d \bar z^{2} + dz^{3} \wedge d \bar z^{3}) .
\end{gather*}
\end{subequations}
The stabiliser of $\Omega_{0}$ in $GL_{+}(\bbr^{6})$ is
$SL(\bbc^{3})$. Therefore each $\Omega \in \Lambda^{3}_{+}V^{*}$ defines a
complex structure $J$ on~$V$, and there is a unique
$\hat \Omega \in \Lambda^{3}_{+}V^{*}$ such that $\Omega + i \hat \Omega$
is a $(3,0)$-form. We say that $\Omega$~defines an \slstr{} on~$V$.
(Reversing the
orientation of $V$ changes the sign of both $J$ and $\hat \Omega$.)
By dimension-counting $\Lambda^{3}_{+}V^{*}$ is an open
subset of $\Lambda^{3}V^{*}$.

Similarly the stabiliser in $GL(\bbr^{6})$ of the pair
$(\Omega_{0}, \omega_{0})$
is $SU(3)$, so we can define an $SU(3)$-structure on $V$ to be a pair
$(\Omega, \omega) \in \Lambda^{3}V^{*} \times \Lambda^{2}V^{*}$ which is
equivalent to $(\Omega_{0}, \omega_{0})$ under some oriented
linear isomorphism
$V \cong \bbr^{6}$. An $SU(3)$-structure naturally defines a complex structure
$J$ as above, and also an inner product
$g(\cdot,\cdot) = \omega(\cdot, J\cdot)$.
With respect to the Hodge star defined by this metric $\hat\Omega = *\Omega$.

\begin{defn}
\label{cy3def}
An \emph{$SU(3)$-structure} on an oriented manifold $X^{6}$ is a pair of forms
$(\Omega, \omega) \in \Omega^{3}(X) \times \Omega^{2}(X)$
which defines an $SU(3)$-structure on each tangent space.

$(\Omega, \omega)$ is said to be a \emph{\cystr} if
$\nabla\Omega = 0, \nabla\omega = 0$ with respect to the metric induced by
$(\Omega,\omega)$. $X$ equipped with the structure $(\Omega,\omega)$ and
the associated Riemannian metric is called a \emph{Calabi-Yau 3-fold}.
\end{defn}

If $X$ is a Calabi-Yau $3$-fold in this sense then by proposition
\ref{simpleholprop} the holonomy of the induced metric is contained in $SU(3)$,
and conversely any metric with holonomy contained in $SU(3)$ is induced
by some \cystr. The almost complex structure
$J$ defined by a \cystr{} is integrable, and the metric is Kähler. Moreover
$\Omega + i\hat\Omega$ is a global holomorphic $(3,0)$-form,
and the metric is Ricci-flat.

\cystr s on $X^{6}$ are equivalent to torsion-free cylindrical \gstr s on
$X \times \bbr$ in the following sense.

\begin{defn}
\label{exactcylg2def}
Let $X^{6}$ a compact oriented manifold,
and denote by $t$ the $\bbr$-coordinate on the cylinder $X \times \bbr$.
A \gstr{} $\varphi$ on $X \times \bbr$ is \emph{cylindrical} if it is
translation-invariant and the associated metric is a product metric
$g_{\varphi} = g_{X} + dt^{2}$, for some metric $g_{X}$ on $X$.
\end{defn}

Comparing the point-wise models (\ref{g2formeq}) and (\ref{su3formseq})
it is easy to see that
$(\Omega, \omega)$ is an $SU(3)$-structure on $X$ with metric $g_{X}$
if and only if the translation-invariant
\gstr{} $\varphi = \Omega + dt \wedge \omega$ on $X \times \bbr$
defines the product metric $g_{X} + dt^{2}$.
$Hol(g_{X} + dt^{2}) \subseteq G_{2}$ if and only if 
$Hol(g_{X}) \subseteq SU(3)$, so

\begin{prop}
\label{cylg2prop}
$(\Omega,\omega)$ is a \cystr{} on $X^{6}$ if and only if
$\Omega + dt \wedge \omega$ is a torsion-free cylindrical \gstr{} on
$X \times \bbr$.
\end{prop}

\begin{rmk}
\label{cyldualrmk}
If $\varphi = \Omega + dt \wedge \omega$ is a cylindrical \gstr{} then
\begin{equation*}
*_{\varphi}\varphi = \half \omega^{2} - dt \wedge \hat\Omega .
\end{equation*}
\end{rmk}

\subsection{Manifolds with cylindrical ends}
\label{cylendssub}

We define manifolds with cylindrical ends and their long exact sequence
for cohomology relative to the boundary. We also define
exponentially asymptotically cylindrical (EAC) metrics and \gstr s.

\begin{defn}
\label{cylenddef}
A manifold $M$ is said to have \emph{cylindrical ends} if $M$ is written as
union of two pieces $M_{0}$ and $M_{\infty}$ with common boundary $X$, where
$M_{0}$ is compact, and $M_{\infty}$ is identified with $X \times \bbrp$
by a diffeomorphism (identifying $\partial M_{\infty}$ with $X \times \{0\}$).
$X$~is called the \emph{cross-section} of~$M$.

A \emph{cylindrical coordinate} on $M$ is a smooth function $t : M \to \bbr$
which is equal to the $\bbrp$-coordinate on
$M_{\infty}$ and is negative in the interior of~$M_{0}$.
\end{defn}

The interior of any compact manifold with boundary can be
considered as a manifold with cylindrical ends
by the collar neighbourhood theorem. Conversely
if $M$ has cylindrical ends then we can compactify $M$ by including it
in $\overline M = M_{0} \cup (X \times [0,\infty])$, i.e.
by `adding a copy of $X$ at infinity'.
The cohomology of $\overline M$ relative to its boundary can be identified
with $H^{*}_{cpt}(M)$, the cohomology of the complex $\Omega^{*}_{cpt}(M)$
of compactly supported forms.
The long exact sequence for relative cohomology of $\overline M$
can be written as
\begin{equation}
\label{relexacteq}
\cdots \longrightarrow H^{m-1}(X) \stackrel{\partial}{\longrightarrow}
H^{m}_{cpt}(M) \stackrel{\inc}{\longrightarrow} H^{m}(M)
\stackrel{j^{*}}{\longrightarrow} H^{m}(X) \longrightarrow \cdots
\end{equation}
$\inc : H^{m}_{cpt}(M) \to H^{m}(M)$ is induced
by the inclusion $\Omega^{*}_{cpt}(M) \hookrightarrow \Omega^{*}(M)$.
The image of $e$ is the subspace of cohomology classes with compact
representatives.

\begin{defn}
Let
$H^{m}_{0}(M) = \im \left(\inc : H^{m}_{cpt}(M) \to H^{m}(M) \right)$.
\end{defn}

Cylindrical ends allow us to define a notion of asymptotic
translation-invariance.

\begin{defn}
A tensor field 
on $X \times \bbr$ is called
\emph{translation-invariant} if it is invariant under the obvious
$\bbr$-action on $X \times \bbr$.
\end{defn}

\begin{defn}
Let $M$ be a manifold with cylindrical ends. Call a smooth function
$\rho : M \to \bbr$ a \emph{cut-off function for the cylinder}
if it is $0$ on the compact piece $M_{0}$ and $1$ outside a compact subset
of $M$.
\end{defn}

If $s_{\infty}$ is a section of a vector bundle associated to the tangent
bundle on the cylinder $X \times \bbr$ and $\rho$ is a cut-off function for
the cylinder on $M$ then $\rho s_{\infty}$ can be considered
to be a section of the corresponding vector bundle over $M$.

\begin{defn}
\label{asytensordef}
Let $M$ be a manifold with cylindrical ends and cross-section $X$.
Pick an arbitrary product metric $g_{X} + dt^{2}$ on $X \times \bbr$, and
a cut-off function $\rho$ for the cylinder.
A section $s$ of a vector bundle associated to $TM$ is said to
be \emph{decaying} if $\norm{\nabla^{k}s} \to 0$ uniformly on $X$
as $t \to \infty$ for all $k \geq 0$.
$s$ is said to be \emph{asymptotic} to a translation-invariant section
$s_{\infty}$ of the corresponding bundle on $X \times \bbr$ if
$s - \rho s_{\infty}$ decays.

\label{expasytensordef}
Similarly $s$ is said to be \emph{exponentially decaying} with
rate $\delta > 0$ if  $e^{\delta t} \norm{\nabla^{k}s}$ is bounded
on $M_{\infty}$ for all $k \geq 0$, and \emph{exponentially asymptotic}
to a translation-invariant section $s_{\infty}$ if
$s - \rho s_{\infty}$ decays exponentially.
Denote by $C^{\infty}_{\delta}(E)$ the space of sections of $E$
which decay exponentially with rate~$\delta$.
\end{defn}

The natural choice of topology on $C^{\infty}_{\delta}(E)$ is to
require the linear isomorphism
$C^{\infty}_{\delta}(E) \to C^{\infty}(E)$, $s \mapsto e^{\delta t}s$
to be a homeomorphism.

\begin{defn}
A metric $g$ on a manifold $M$ with cylindrical ends is said to be \emph{EAC}
if it is exponentially asymptotic to a product metric
$g_{X} + dt^{2}$ on $X \times \bbrp$.
An EAC manifold is a manifold with cylindrical ends
equipped with an EAC metric.
\end{defn}

\begin{defn}
\label{cyldiffdef}
Let $M$ be a manifold with cylindrical ends and cross-section $X$.
A diffeomorphism $\Psi_{\infty}$ of the cylinder $X \times \bbr$ is said to be
cylindrical if it is of the form
\[ \Psi_{\infty}(x,t) = (\Xi(x), t + h) , \]
where $\Xi$ is a diffeomorphism of $X$ and $h \in \bbr$.
A diffeomorphism $\Psi$ of $M$ is said to be \emph{EAC} with rate $\delta > 0$
if there is a cylindrical diffeomorphism
$\Psi_{\infty}$ of $X \times \bbr$, a real $T > 0$ and an exponentially
decaying vector field $V$ on $M$ such that on $X \times (T, \infty)$
\[ \Psi = (\exp V) \circ \Psi_{\infty} . \]
\end{defn}

\begin{defn}
\label{cylg2def}
Let $M^{7}$ be a connected oriented manifold with cylindrical ends
and cross-section $X^{6}$.
A \gstr{} $\varphi$ on $M$ is said to be
EAC if it is exponentially asymptotic
to a cylindrical \gstr{} on $X \times \bbr$
(cf. definition~\ref{exactcylg2def}).
$M$~equipped with a torsion-free EAC
\gstr{} and the associated metric is called an
\emph{EAC \gmfd}.
\end{defn}

If $\varphi$ is a torsion-free EAC \gstr{} then note that the associated
metric $g_{\varphi}$ is EAC, and that by proposition
\ref{cylg2prop} the asymptotic limit defines a Calabi-Yau structure on the
cross-section~$X$.

The next theorem implies that an EAC \gmfd{} is not interesting unless it has
a single end. The theorem can be proved using the Cheeger-Gromoll splitting
theorem, and is also proved using more elementary methods by Salur
\cite{salur05}.

\begin{thm}
\label{salurthm}
Let $M$ be an orientable connected asymptotically 
cylindrical Ricci-flat manifold. Then either
$M$ has a single end, i.e. its cross-section $X$ is connected,
or $M$ is a cylinder $X \times \bbr$ with a product metric.
\end{thm}

\section{Statement of results}
\label{statesec}

Let $M^{7}$ be a connected oriented manifold with cylindrical ends
and cross-section $X^{6}$.
For $\delta > 0$ let $\cald{X}$ be the space of smooth
torsion-free exponentially
asymptotically cylindrical (EAC) \mbox{\gstr s} with rate $\delta$ on $M$
(see definition \ref{cylg2def}).
$\cald{X}$ has the topology of a subspace of the space of smooth
exponentially asymptotically translation-invariant $3$-forms.

Let $\calp{X} = \bigcup_{\delta>0} \cald{X}$. If $\delta_{1} > \delta_{2} > 0$
then the inclusion
$\mathcal{X}_{\delta_{1}} \hookrightarrow \mathcal{X}_{\delta_{2}}$
is continuous, so we can give $\calp{X}$ the direct limit topology,
i.e. $U \subseteq \calp{X}$ is open if and only if $U \cap \cald{X}$
is open in $\cald{X}$ for all $\delta > 0$. Similarly let $\calp{D}$ be the
group of 
EAC diffeomorphisms of~$M$ with any positive rate
(in the sense of definition \ref{cyldiffdef})
that are isotopic to the identity. $\calp{D}$ acts on $\calp{X}$ by pull-backs,
and the
\emph{moduli space of torsion-free 
EAC \gstr s} on $M$ is the quotient $\calp{M} = \calp{X}/\calp{D}$.

\begin{rmk}
\label{normalisermk}
The definition of an EAC \gstr{} $\varphi$ that is used involves a
normalisation -- if $t$ is the cylindrical coordinate on $M$ then
$\norm{\contrat} \to 1$ uniformly on $X$ as $t \to \infty$
(in the metric defined by $\varphi$),
so a scalar multiple $\lambda \varphi$ is not an EAC \gstr{}.
This normalisation is the most convenient to work with, but a different
choice of normalisation (e.g. that $Vol(X) = 1$ in the induced metric
on the boundary) would of course give the same results.
Another interpretation is that $\bbrp$ acts on the moduli space of unnormalised
EAC \gstr s by rescaling, and that $\calp{M}$ is the resulting quotient.
\end{rmk}

In the compact case theorem \ref{maincptthm} gives a description of the moduli
space of
torsion-free \gstr s using the natural projection map to the de Rham cohomology.
In the EAC case, however, it is not enough to consider
\[ \pi_{H} : \calp{M} \to H^{3}(M), \;\: \varphi\calp{D} \mapsto [\varphi] . \]
We also need to consider the boundary values of $\varphi$ to get an adequate
description.
Any $\varphi \in \calp{X}$ is asymptotic to some $\Omega + dt \wedge \omega$
with $(\Omega, \omega) \in \Omega^{3}(X) \times \Omega^{2}(X)$.
Let
\begin{equation}
\label{derhamprojeq}
\pi_{\mathcal{M}} : \calp{M} \to H^{3}(M) \times H^{2}(X), \;\:
\varphi\calp{D} \mapsto ([\varphi], [\omega]) .
\end{equation}
The main theorem we shall prove is

\begin{thm}
\label{maincylthm}
$\calp{M}$ is a smooth manifold, and
$\pi_{\mathcal{M}} : \calp{M} \to H^{3}(M) \times H^{2}(X)$ is an immersion.
\end{thm}

In order to prove theorem \ref{maincylthm} we will need to understand
the deformations of the `boundary' of an EAC \gmfd, i.e. the deformations of
torsion-free cylindrical \gstr s. By proposition \ref{cylg2prop} this
corresponds to deformations of \cystr s.
Let $\mathcal{Y}$ be the set of \cystr s $(\Omega, \omega)$
on $X$, and $\calx{D}$ the group of diffeomorphisms of $X$ isotopic
to the identity. The \emph{moduli space of Calabi-Yau structures on $X$} is
$\mathcal{N} = \mathcal{Y}/\calx{D}$, and there is a natural
projection to the de Rham cohomology
\begin{equation}
\pi_{\mathcal{N}} : \mathcal{N} \to H^{3}(X) \times H^{2}(X), \;\:
(\Omega,\omega)\calx{D} \mapsto ([\Omega], [\omega]) .
\end{equation}

\begin{thm}
\label{cymodthm}
Let $X^{6}$ be a compact connected Calabi-Yau $3$-fold.
The moduli space $\mathcal{N}$ of \cystr s on $X$ is a manifold,
\begin{equation}
\dim \mathcal{N} = b^{3}(X)+b^{2}(X)-b^{1}(X)-1 ,
\end{equation}
and $\pi_{\mathcal{N}} : \mathcal{N} \to H^{3}(X) \times H^{2}(X)$
is an immersion.
\end{thm}

\begin{rmk}
The definition of a Calabi-Yau 3-fold $X^{6}$ used here allows
$Hol(X)$ to be a proper subgroup of $SU(3)$. If $Hol(X)$ is exactly $SU(3)$
(so $X$ is irreducible as a Riemannian manifold)
then $b^{1}(X) = 0$, and the formula for the dimension simplifies to
$b^{3}(X)+b^{2}(X)-1$.

If $X$ is an irreducible Calabi-Yau manifold then for any Calabi-Yau structure
$(\Omega,\omega)$ on $X$ and $\lambda \in \bbrp$ we can define a torsion-free
product \gstr{} $\varphi = \Omega + \lambda d\theta \wedge \omega$ on
$X \times S^{1}$. The metric defined by $\varphi$ is the product of the
Calabi-Yau metric on $X$ and the metric on $S^{1}$ with radius $\lambda$
(cf. proposition \ref{cylg2prop}). The moduli space
of such torsion-free product \gstr s has dimension
\[ \dim \mathcal{N} + 1 = b^{3}(X) + b^{2}(X) = b^{3}(X \times S^{1}) , \]
which equals the dimension of the moduli space of torsion-free \gstr s
on $X \times S^{1}$ by theorem \ref{maincptthm}.
\end{rmk}

Theorem \ref{cymodthm} is actually a special case of a known result.
The deformation theory for complex manifolds was developed by
Kodaira and Spencer. Tian \cite{tian87} and Todorov \cite{todorov89}
showed independently that on a compact connected Calabi-Yau manifold $X^{2n}$
with holonomy exactly $SU(n)$ the deformations of the complex structure are
`unobstructed'.
This implies that the moduli space of complex structures on $X$ is a manifold
of dimension $2h^{1,n-1}(X)$ ($h^{p,q}(X)$ denote the Hodge numbers of $X$,
i.e. the dimension of the Dolbeault cohomology $H^{p,q}(X)$).
It is easy to deduce from this and Yau's solution of the Calabi conjecture
\cite{yau78} that the moduli space of Calabi-Yau structures on a complex
$n$-fold (in the sense of an integrable complex structure with a Kähler metric
and a holomorphic $(n,0)$-form of norm $1$) is a manifold of dimension
\begin{equation}
\label{tianmodeq}
2h^{1,n-1}(X) + h^{1,1}(X) + h^{n,0}(X)
\end{equation}
(cf. \cite[Section 6.8]{joyce00}). By \cite[Proposition 6.2.6]{joyce00}
$h^{m,0}(X) = 0$ for $0 < m < n$ and $h^{n,0}(X) = 1$ when
$X^{2n}$ is compact connected and $Hol(X) = SU(n)$,
so if $n = 3$ then $b^{3}(X) = 2h^{1,2}(X) + 2$ and
$b^{2}(X) = h^{1,1}(X)$. Thus the expression (\ref{tianmodeq}) for the
dimension of the moduli space can be rewritten as $b^{3}(X) + b^{2}(X) - 1$
when $n=3$, which agrees with the formula stated in theorem~\ref{cymodthm}.

We will give a different proof of theorem \ref{cymodthm} in section
\ref{cymodsec}. We produce pre-moduli spaces by an elementary application
of the implicit function theorem, extending arguments of Hitchin in
\cite{hitchin00}. These pre-moduli spaces are also used in the proof
of theorem \ref{maincylthm}.

In subsection \ref{propertiessub} we look at some local properties
of $\calp{M}$. Its dimension is given by

\begin{prop}
\label{dimmodprop}
$\dim \calp{M} = b^{4}(M) + \half b^{3}(X) - b^{1}(M) - 1$.
\end{prop}

We also study the properties of the boundary map on $\calp{M}$,
i.e. the map $B : \calp{M} \to \mathcal{N}$
which sends a \gstr{} on $M$ to the \cystr{} on $X$ defined by its asymptotic
limit. Denote by $A^{m} \subseteq H^{m}(X)$ the image of the
pull-back map $j^{*} : H^{m}(M) \to H^{m}(X)$ in the long exact sequence
for relative cohomology (\ref{relexacteq}).
If $\varphi$ is asymptotic to $\Omega + dt \wedge \omega$ then
\begin{gather*}
[\Omega] = j^{*}([\varphi]) \in A^{3} , \\
\half [\omega^{2}] = j^{*}([*_{\varphi}\varphi]) \in A^{4} ,
\end{gather*}
so the image of $B : \calp{M} \to \mathcal{N}$ is contained in
\begin{equation}
\label{nadefeq}
\mathcal{N}_{A} = \{ (\Omega, \omega)\calx{D} \in \mathcal{N} :
[\Omega] \in A^{3}, [\omega^{2}] \in A^{4} \} .
\end{equation}
It turns out that -- locally at least -- these necessary conditions for
a point to be in the image are also sufficient.

\begin{thm}
\label{modsubmersethm}
The image of
\begin{equation}
\label{modsubmerseeq}
B : \calp{M} \to \mathcal{N}_{A}
\end{equation}
is open in $\mathcal{N}_{A}$ and a submanifold of $\mathcal{N}$.
The map is a submersion onto its image.
\end{thm}

Since the methods used are entirely local they do not tell us anything
about the global properties of $\calp{M}$ or the image of
(\ref{modsubmerseeq}).

We will show that the fibres of the submersion (\ref{modsubmerseeq})
are locally diffeomorphic to the compactly supported subspace
$H^{3}_{0}(M) \subseteq H^{3}(M)$. The fibre over
$(\Omega,\omega)$ corresponds to the moduli space of torsion-free \gstr s
asymptotic to $\Omega + dt \wedge \omega$. Thus

\begin{cor}
\label{fixboundarycor}
The moduli space of torsion-free \gstr s on $M$ exponentially asymptotic to a
fixed cylindrical \gstr{} on $X \times \bbr$ is a manifold locally
diffeomorphic to $H^{3}_{0}(M)$.
\end{cor}

In the proof of proposition \ref{dimmodprop} we find that
$\dim H^{3}_{0}(M) = b^{3}(M) - \half b^{3}(X)$, so
\[ \dim \mathcal{N}_{A} = b^{4}(M) - b^{3}(M) + b^{3}(X) - b^{1}(M) - 1 . \]

Theorem \ref{salurthm} implies that if $M$ is a \gmfd{} then either $M$ is
a cylinder $X \times \bbr$ (with a product metric) or $M$ has a single end.
If $M$ is a cylinder $X \times \bbr$ then the only possible torsion-free
\gstr{} asymptotic to
a given cylindrical \gstr{} $\varphi_{\infty}$ is $\varphi_{\infty}$ itself, so
the moduli space of asymptotically cylindrical torsion-free \gstr s on $M$ is
equivalent to the moduli space of Calabi-Yau structures on $X$ (we can compute
that $H^{m}_{0}(X \times \bbr) = 0$ for all $m$, so this agrees with
corollary \ref{fixboundarycor}).
The moduli space will therefore only be interesting when $M$ has a single end.
We will not need to assume this in the proof of theorem \ref{maincylthm}
though.

Finally, in section \ref{topcondsec} we find a topological condition for
when the holonomy group
of an EAC \gmfd{} is exactly $G_{2}$, as opposed to a proper subgroup.
For compact \gmfd s it is well-known that the holonomy is exactly
$G_{2}$ if and only if the fundamental group is finite. In the EAC case
we need to take into account that a product cylinder $X^{6} \times \bbr$ may
have finite fundamental group, but cannot have holonomy $G_{2}$.
The correct statement is

\begin{thm}
\label{maintopcondthm}
Let $M^{7}$ be an EAC \gmfd. Then $Hol(M) = G_{2}$ if and only if
the fundamental group $\pi_{1}(M)$ is finite and neither
$M$ nor any double cover of $M$ is homeomorphic to a cylinder.
\end{thm}

\section{Deformations of compact Calabi-Yau $3$-folds}
\label{cymodsec}

In \cite{hitchin00} Hitchin uses elementary methods to construct the moduli
space of torsion-free \slstr s on a compact manifold $X^{6}$.
In a sense this 
provides `half' the deformation theory for compact Calabi-Yau 3-folds.
In this section we show how to extend Hitchin's arguments to construct
pre-moduli spaces of \cystr s, which we require for the proof of the
main theorem \ref{maincylthm}.
This also gives an elementary proof of
theorem \ref{cymodthm}, stating that the moduli space of \cystr s
on $X$ is a manifold.

\subsection{Harmonic forms and holonomy}
\label{harmholsub}

We first review how a holonomy constraint on a Riemannian manifold gives rise
to decompositions of harmonic forms, similar to the Hodge decomposition on
a Kähler manifold.
This is explained in more detail in \cite[Section $3.5$]{joyce00}.

Let $H$ be a closed subgroup of $SO(n)$, and $M^{n}$ an oriented Riemannian
manifold with $Hol(M) \subseteq H$, equipped with a corresponding $H$-structure.
Suppose that $\Lambda^{m}\bbr^{n}$ decomposes as an orthogonal direct sum of
subrepresentations $\Lambda^{m}\bbr^{n} = \bigoplus \Lambda^{m}_{d}\bbr^{n}$
under the action of $H$ (we will indicate the rank of the subrepresentations
by the index $d$). Then there is a
corresponding $H$-invariant decomposition of the exterior product bundle
\begin{equation}
\label{extdecompeq}
\Lambda^{m}T^{*}M = \bigoplus \Lambda^{m}_{d}T^{*}M .
\end{equation}
We write $\Omega^{m}_{d}(M)$ for the space sections of
$\Lambda^{m}_{d}T^{*}M$, the `forms of type~d'. We define projections
$\pi_{d} : \Lambda^{m}T^{*}M \to \Lambda^{m}_{d}T^{*}M$. These induce 
maps $\pi_{d} : \Omega^{m}(M) \to \Omega^{m}_{d}(M)$ on the sections, and
allow us to decompose forms into type components.
As observed by Chern in \cite{chern57}, the condition that $Hol(M) \subseteq H$
ensures that the Hodge Laplacian respects these type decompositions.

\begin{prop}[{\cite[Theorem $3.5.3$]{joyce00}}]
\label{hodgelapprop}
Let $M^{n}$ be a Riemannian manifold with $Hol(M) \subseteq H$.
If $\Lambda^{m}_{d}T^{*}M$ is an $H$-invariant subbundle of
$\Lambda^{m}T^{*}M$ then $\triangle$ commutes with $\pi_{d}$ on
$\Omega^{m}(M)$, and maps $\Omega^{m}_{d}(M)$ to itself.
Moreover, if $\Lambda^{k}_{e}T^{*}M$ is an $H$-invariant subbundle of
$\Lambda^{k}T^{*}M$ and
$\phi : \Lambda^{m}_{d}T^{*}M \to \Lambda^{k}_{e}T^{*}M$
is $H$-equivariant then the diagram below commutes. 
\begin{diagram}[PostScript=Rikicki]
\Omega^{m}_{d}(M) &\rTo^{\phi} & \Omega^{k}_{e}(M) \\
\dTo^{\triangle} & & \dTo_{\triangle} \\
\Omega^{m}_{d}(M) &\rTo^{\phi} & \Omega^{k}_{e}(M) 
\end{diagram}
\end{prop}

It follows that given an $H$-invariant decomposition (\ref{extdecompeq})
of $\Lambda^{m}T^{*}M$ into subbundles there is a corresponding
decomposition of the harmonic forms
\[ \harm^{m} = \bigoplus \harm^{m}_{d} . \]
If $M$ is compact then by Hodge theory the natural map $\harm^{m} \to H^{m}(M)$
is an isomorphism. If we let $H^{m}_{d}(M)$ be the image of $\harm^{m}_{d}$
then we obtain a decomposition of the de Rham cohomology
\begin{equation}
\label{hmdecompeq}
H^{m}(M) = \bigoplus H^{m}_{d}(M) .
\end{equation}

In this section we consider a Calabi-Yau 3-fold $X^{6}$.
The standard representation of $SU(3)$ on $\bbr^{6}$ is irreducible,
and $\Lambda^{m}\bbr^{6}$ decomposes~as
\begin{align*}
\Lambda^{2}\bbr^{6} &= \Lambda^{2}_{1}\bbr^{6}
\oplus \Lambda^{2}_{6}\bbr^{6} \oplus \Lambda^{2}_{8}\bbr^{6} , \\
\Lambda^{3}\bbr^{6} &= \Lambda^{3}_{1\oplus1}\bbr^{6} 
\oplus \Lambda^{3}_{6}\bbr^{6} \oplus \Lambda^{3}_{12}\bbr^{6} , \\
\Lambda^{4}\bbr^{6} &= \Lambda^{4}_{1}\bbr^{6}
\oplus \Lambda^{4}_{6}\bbr^{6} \oplus \Lambda^{4}_{8}\bbr^{6} .
\end{align*}
Each of the subrepresentations $\Lambda^{m}_{d}\bbr^{6}$ is irreducible,
but $\Lambda^{3}_{1\oplus1}\bbr^{6}$ is trivial of rank $2$.
The corresponding decompositions of the
exterior cotangent bundles of $X$ are related to the Hodge decomposition, e.g.
$\Lambda^{2}_{1}T^{*}X \oplus \Lambda^{2}_{8}T^{*}X$
is the bundle of real $(1,1)$-forms, while
$\Lambda^{2}_{6}T^{*}X$ consists of the real and imaginary parts of forms of
type $(2,0)$.

\subsection{Pre-moduli space of Calabi-Yau structures}
\label{cy3premodsub}

Let $X^{6}$ a compact connected oriented manifold.
Recall that in subsection \ref{cydefsub} we defined a Calabi-Yau structure on $X$
in terms of a pair of forms
$(\Omega,\omega) \in \Omega^{3}(X) \times \Omega^{2}(X)$.
As defined in section \ref{statesec} the
\emph{moduli space of Calabi-Yau structures}
on $X$ is $\mathcal{N} = \mathcal{Y}/\calx{D}$, where
$\mathcal{Y}$ is the set of \cystr s $(\Omega, \omega)$
on $X^{6}$, and $\calx{D}$ is the group of
diffeomorphisms of $X$ isotopic to the identity.
To prove theorem \ref{cymodthm} we find pre-moduli spaces of \cystr s,
i.e. manifolds $\mathcal{Q} \subseteq \mathcal{Y}$ that are homeomorphic
to open sets in $\mathcal{N}$, and can therefore be used as charts.

On $6$-dimensional manifolds we have the following convenient characterisation
of \cystr s (cf. Hitchin \cite[Section 2]{hitchin97})

\begin{lem}
\label{cy3charlem}
Let $X^{6}$ an oriented manifold, and
$(\Omega,\omega) \in \Omega^{3}(X) \times \Omega^{2}(X)$. Suppose that
$\Omega$ is stable, so that it defines an almost complex structure $J$ and
a $3$-form $\hat\Omega$. Then the
following conditions are sufficient to ensure that $(\Omega,\omega)$ is
a \cystr :
\begin{enumerate}
\item \label{volitem} $\quart \Omega \wedge \hat\Omega = \sixth \omega^{3}$
\item \label{type11item} $\Omega \wedge \omega = 0$
\item \label{posdefitem}
$g(\cdot,\cdot) = \omega(\cdot, J\cdot)$ is positive-definite
\item \label{hatitem} $d\Omega = d\hat\Omega = 0$
\item \label{kahleritem} $d\omega = 0$
\end{enumerate}
\end{lem}

Recall that $\hat\Omega$ is the unique $3$-form such that
$\Omega + i\hat\Omega$ has type $(3,0)$ with respect to the complex structure
defined by $\Omega$. $\hat\Omega$ depends smoothly on the stable form $\Omega$.
The derivative of $\Omega \mapsto \Omega \wedge \hat\Omega$ must be
proportional to $\cdot \wedge \hat\Omega$ since it is $SU(3)$-equivariant.
The proportionality constant is $2$
as $\hat\Omega$ is homogeneous of degree $1$ in $\Omega$
(cf. \cite[page 10]{hitchin00}).

Now consider a fixed \cystr{} $(\Omega,\omega)$.
Any tangent $(\sigma,\tau)$ to a path of
$SU(3)$-structures through $(\Omega,\omega)$ must satisfy the linearisation
of the point-wise algebraic conditions \ref{volitem} and \ref{type11item}
in lemma \ref{cy3charlem}, i.e.
\begin{subequations}
\label{leqs}
\begin{gather}
\label{loneeq}
L_{1}(\sigma,\tau) = \sigma \wedge \hat \Omega - \tau \wedge \omega^{2} \;= 0,\\
\label{ltwoeq}
L_{2}(\sigma,\tau) = \sigma \wedge \omega + \Omega \wedge \tau = 0 .
\end{gather}
\end{subequations}

\begin{defn}
\label{hsudef}
Let
\[ \harm_{SU} = \{ (\sigma,\tau) \in \harm^{3} \times \harm^{2} :
L_{1}(\sigma,\tau) = L_{2}(\sigma,\tau) = 0 \} . \]
\end{defn}

The natural map
$\pi_{\mathcal{N}} : \harm_{SU} \to H^{3}(X) \times H^{2}(X), \;
(\sigma,\tau) \mapsto ([\sigma], [\tau])$
is injective by Hodge theory for compact manifolds. Proposition
\ref{hodgelapprop} implies that
$L_{1} : \harm^{3} \times \harm^{2} \to \harm^{6}$ and
$L_{2} : \harm^{3} \times \harm^{2} \to \harm^{5}$ are surjective, so
\begin{equation}
\dim \harm_{SU} = b^{3}(X) + b^{2}(X) - b^{1}(X) - 1 .
\end{equation}


Theorem \ref{cymodthm} will follow from the existence of a pre-moduli space near
each $(\Omega,\omega)$ with tangent space $\harm_{SU}$.

\begin{prop}
\label{cypremodprop}
For any $(\Omega, \omega) \in \mathcal{Y}$ there is a manifold
$\mathcal{Q} \subseteq \mathcal{Y}$ near $(\Omega,\omega)$
such that the natural map
$\mathcal{Q} \to \mathcal{N}$ is a homeomorphism onto an open subset.
The tangent space of $\mathcal{Q}$ at $(\Omega,\omega)$ is $\harm_{SU}$.
\end{prop}

We find $\mathcal{Q}$ with the desired properties using a slice
construction.
Pick some $k \geq 1$, $\alpha \in (0,1)$ and let $\ztre_{k} \times \ztwo_{k}$
be the space of pairs of closed Hölder $\holda{k}$ $3$- and $2$-forms.
$\mathcal{Y} \hookrightarrow \ztre_{k} \times \ztwo_{k}$ continuously.
We find a direct complement $\slt$ in $\ztre_{k} \times \ztwo_{k}$ to the
tangent space of the $\calx{D}$-orbit at $(\Omega,\omega)$,
and use a small neighbourhood $\mathcal{S}$ of $(\Omega,\omega)$ in the affine
space $(\Omega,\omega) + \slt$ as a \emph{slice} for the $\calx{D}$-action.

Let $\mathcal{Q} \subseteq \mathcal{S}$ be the subspace of elements which
define \cystr s.
Subsection \ref{slsub} summarises Hitchin's construction of the moduli
space of torsion-free \slstr s (stable forms $\Omega$ satisfying condition
\ref{hatitem} in lemma \ref{cy3charlem}). In subsection
\ref{cypremodproofsub} we extend those arguments in order to prove that
$\mathcal{Q} \subseteq \mathcal{S}$ is a submanifold by
an application of the implicit function theorem.

By elliptic regularity the elements of $\mathcal{Q}$ are smooth,
and slice arguments show that $\mathcal{Q}$ is homeomorphic to a neighbourhood
of $(\Omega,\omega)\calx{D}$ in~$\mathcal{N}$.
Such arguments will be explained in detail in the more complicated EAC case
(cf. propositions \ref{cylregprop} and \ref{qsliceprop}).
Thus $\mathcal{Q}$ satisfies the conditions of proposition \ref{cypremodprop}.

We will use the deformation theory of \cystr s in the proof of the main theorem
\ref{maincylthm} on the moduli space of torsion-free EAC \gstr s.
We will also require the following property of the pre-moduli spaces
of \cystr s.

\begin{prop}
\label{autinvprop}
Let $\mathcal{Q}$ be the pre-moduli space of \cystr s near $(\Omega, \omega)$.
If $\mathcal{Q}$ is taken sufficiently small then all elements of $\mathcal{Q}$
have the same stabiliser in $\calx{D}$.
\end{prop}

\begin{proof}
Let $\mathcal{I} \subseteq \calx{D}$ be the stabiliser of $(\Omega,\omega)$.
$\mathcal{I}$ is contained in the isometry group of a Riemannian metric,
so it is a compact Lie group. By shrinking $\mathcal{Q}$ we may assume that
it is mapped to itself by $\mathcal{I}$. Since
$\pi_{\mathcal{N}} : \mathcal{Q} \to H^{3}(X) \times H^{2}(X)$ is injective
and $\calx{D}$-invariant $\mathcal{I}$ acts trivially on $\mathcal{Q}$.

Conversely, it follows from \cite[Theorem $7.1(2)$]{ebin70} that if
$\phi \in \calx{D}$ fixes an element of $\mathcal{Q}$ sufficiently close
to $(\Omega,\omega)$ then $\phi \in \mathcal{I}$.
\end{proof}

\subsection{Torsion-free $SL(\bbc^{3})$-structures}
\label{slsub}

Recall from page \pageref{sl3formeq} that an \slstr{} on an oriented dimension
$6$ vector space $V$ is defined by a stable $3$-form
$\Omega \in \Lambda^{3}_{+}V^{*}$.

\begin{defn}
An \emph{\slstr} on an oriented manifold $X^{6}$ is a section $\Omega$ of
$\Lambda^{3}_{+}T^{*}X$.
$\Omega$ is \emph{torsion-free} if $d\Omega = d\hat\Omega = 0$.
\end{defn}

If $\Omega$ is torsion-free then so is the almost complex structure $J$ it
defines, and $\Omega + i\hat\Omega$ is a global
holomorphic $(3,0)$-form. 
A torsion-free \slstr{} is therefore equivalent to a complex structure
with trivial canonical bundle, together with a choice of trivialisation.

The next two propositions are a summary of Sections $6.1$ and $6.2$
in~\cite{hitchin00}. Let $X^{6}$ be a compact oriented manifold.
Fix a torsion-free \slstr{} $\Omega$ on $X$,
and pick an arbitrary Riemannian metric 
that is Hermitian with respect to the complex structure defined by $\Omega$.
Take $k \geq 1$, $\alpha \in (0,1)$, and let $\ztre_{k}$ be the space of closed
$\holda{k}$ $3$-forms.
Abbreviate $\Lambda^{m}T^{*}X$ to $\Lambda^{m}$.
The complex structure $J$ defines a vector bundle splitting
\begin{equation*}
\Lambda^{2} = \Lambda^{2}_{6} \oplus \Lambda^{2}_{9} ,
\end{equation*}
where $\Lambda^{2}_{9}$ denotes the bundle of real $(1,1)$-forms,
and $\Lambda^{2}_{6}$ has type $(2,0)+(0,2)$.

\begin{prop}
\label{slsliceprop}
There is an \ltwoorth{} direct sum decomposition
\begin{equation*}
\ztre_{k} = \harm^{3} \oplus W_{1} \oplus d\holda{k+1}(\Lambda^{2}_{6}) ,
\end{equation*}
and the projections are bounded in the Hölder $\holda{k}$-norm.
\end{prop}

Let $P_{1} : \holda{k}(\Lambda^{3}) \to W_{1}$ be the \ltwoorth{}
projection. If
$\beta$ is sufficiently close to $\Omega$ then $\beta$ is stable,
and $\hat\beta$ is well-defined.
On a neighbourhood of $\Omega$ in $\ztre_{k}$ we 
define
\begin{equation*}
\label{fonedefeq}
F_{1}(\beta) = P_{1}(*\hat\beta) .
\end{equation*}

\begin{prop}
\label{foneprop}
The derivative $(DF_{1})_{\Omega} : \ztre_{k} \to W_{1}$ is
$0$ on $\harm^{3} \oplus d\holda{k+1}(\Lambda^{2}_{6})$ and bijective
on $W_{1}$. Furthermore for $\beta \in \ztre_{k}$ sufficiently close
to~$\Omega$
\begin{equation*}
F_{1}(\beta) = 0 \Leftrightarrow d\hat\beta = 0 .
\end{equation*}
\end{prop}

In \cite{hitchin00} the content of the above two propositions is used with a
slice argument to 
construct a moduli space of torsion-free \slstr s.
In the next subsection we
extend the argument to prove proposition \ref{cypremodprop}.

\subsection{Proof of proposition \ref{cypremodprop}}
\label{cypremodproofsub}

We now explain how to define a slice for the $\calx{D}$-action at
$(\Omega,\omega) \in \mathcal{Y}$. We find a function with surjective
derivative whose zero set in the slice is precisely the subspace of
\cystr s $\mathcal{Q}$. Together with elliptic regularity and slice arguments
this proves proposition \ref{cypremodprop}, and hence theorem \ref{cymodthm}.

Abbreviate $\Lambda^{m}T^{*}X$ to $\Lambda^{m}$. As described in subsection
\ref{harmholsub} the \cystr{} $(\Omega,\omega)$ induces decompositions
$\Lambda^{2} = \Lambda^{2}_{1} \oplus \Lambda^{2}_{6} \oplus \Lambda^{2}_{8}$
etc.

The tangent space to the $\calx{D}$-orbit at $(\Omega, \omega)$ in
$\ztre_{k} \times \ztwo_{k}$ is
\[ \orbt =
\{ (d(V\lrcorner\Omega), d(V\lrcorner \omega)) : V \in \holda{k+1}(TX) \} . \]
The sections of $\Lambda^{2}_{6}$ are precisely $V \lrcorner \Omega$ for
vector fields $V$, so by proposition \ref{slsliceprop} we may take
$\slt =(\harm^{3} \oplus W_{1}) \times \ztwo_{k}$ as a complement
of $\orbt$ in $\ztre_{k} \times \ztwo_{k}$. It is not clear
\textit{a priori} that $\slt \cap \orbt = 0$, but it will turn out
to be so (corollary \ref{transcor}).
We use a small neighbourhood $\mathcal{S}$ of $(\Omega,\omega)$ in the
affine space $(\Omega,\omega) + \slt$ as a slice for the $\calx{D}$-action.

The pre-moduli space $\mathcal{Q} \subseteq \mathcal{S}$ is the subspace of
elements which define \cystr s.
By lemma \ref{cy3charlem} $\mathcal{Q}$ is the zero set of
\[ \mathcal{S} \to \holda{k-1}(\Lambda^{4}) \times \holda{k}(\Lambda^{5})
\times \holda{k}(\Lambda^{6}), \;\:
(\beta, \gamma) \mapsto (d\hat\beta, \beta \wedge \gamma,
\quart \beta \wedge \hat\beta - \sixth \gamma^{3}) , \]
but this function does not have surjective derivative.
Part of the work of obtaining a more appropriate function
is already done -- we can replace the first component by $F_{1}$ defined
in (\ref{fonedefeq}). We need to find a more suitable second component.

$(\Omega, \omega)$ defines a complex structure on $X$, and in particular
gives us the conjugate differential $d^{c} = i(\bar\partial - \partial)$.
Note that $dd^{c} = -d^{c}d = 2i\partial\bar\partial$.
Since $X$ is a Kähler manifold the $dd^{c}$-lemma holds.

\begin{thm}
Any real exact form of type $(1,1)$ on a compact Kähler manifold
$X$ is $dd^{c}$-exact.
\end{thm}

\begin{prop}
There is an \ltwoorth{} direct sum decomposition
\begin{equation}
\label{zfiveeq}
\zfive_{k} = \harm^{5} \oplus d\holda{k+1}(\Lambda^{4}_{6})
\oplus W_{2},
\end{equation}
where
$W_{2} = \{ d\eta \in d\holda{k+1}(\Lambda^{4}) : \pi_{6}d^{*}d\eta = 0 \}$.
The projection
\begin{equation}
\label{p2defeq}
P_{2} : \holda{k}(\Lambda^{5}) \to
\harm^{5} \oplus d\holda{k+1}(\Lambda^{4}_{6})
\end{equation}
is bounded in Hölder $\holda{k}$-norm.
\end{prop}

\begin{proof}
The operator $\pi_{6}d^{*} \oplus d \oplus d^{c} : \Gamma(\Lambda^{5})
\to \Gamma(\Lambda^{4}_{6}\oplus\Lambda^{6}\oplus\Lambda^{6})$ is
over\-determined elliptic, and its formal adjoint is $d+d^{*}+d^{c*}$.
Hence the operator $d\pi_{6}d^{*} + d^{*}d + d^{c*}d^{c}$ is elliptic 
and formally self-adjoint on sections of $\Lambda^{5}$, so
\begin{multline*}
\holda{k}(\Lambda^{5}) =
(d\pi_{6}d^{*} + d^{*}d + d^{c*}d^{c})\holda{k+2}(\Lambda^{5}) \oplus
\ker (d\pi_{6}d^{*} + d^{*}d + d^{c*}d^{c}) \\
\subseteq d\holda{k+1}(\Lambda^{4}_{6}) + d^{*}\holda{k+1}(\Lambda^{6}) +
d^{c*}\holda{k+1}(\Lambda^{6}) + \ker{\pi_{6}d^{*}} \subseteq
d\holda{k+1}(\Lambda^{4}_{6}) + \ker{\pi_{6}d^{*}} ,
\end{multline*}
and the last sum is clearly direct. $\ker{\pi_{6}d^{*}}$ contains $\harm^{5}$
and $d^{*}\holda{k+1}(\Lambda^{6})$, and hence
splits as $\harm^{5} \oplus W_{2} \oplus d^{*}\holda{k+1}(\Lambda^{6})$.
\end{proof}

\noindent
$d^{c}$ gives a convenient characterisation of the splitting (\ref{zfiveeq}).

\begin{prop}
\label{dccharprop}
$\Omega \wedge \ztwo_{k} = \harm^{5} \oplus d\holda{k+1}(\Lambda^{4}_{6})
= \zfive_{k} \cap \ker d^{c}$.
\end{prop}

\begin{proof}
$\Omega \wedge \harm^{2} = \harm^{5}$ by proposition \ref{hodgelapprop}, while
\[ \Omega \wedge d\holda{k+1}(\Lambda^{1}) =
d(\Omega \wedge \holda{k+1}(\Lambda^{1})) = d\holda{k+1}(\Lambda^{4}_{6}) . \]

A real $4$-form has type $(2,2)$ if and only if its $\Lambda^{4}_{6}$ part
vanishes.
For $d\eta \in d\holda{k+1}(\Lambda^{4})$ applying the $dd^{c}$-lemma
therefore gives
\begin{equation*}
\pi_{6}d^{*}d\eta = 0 \Leftrightarrow
d^{*}d\eta \in d^{*}d^{c*}\holda{k+1}(\Lambda^{6})
\Leftrightarrow d\eta \in P_{E}d^{c*}\holda{k+1}(\Lambda^{6}) ,
\end{equation*}
where $P_{E} : \holda{k}(\Lambda^{5}) \to d\holda{k+1}(\Lambda^{4})$ is the
\ltwoorth{} projection to the exact forms.
Thus $W_{2} = P_{E}d^{c*}\holda{k+1}(\Lambda^{6})$, which is the \ltwoorth{}
complement
to $d\holda{k+1}(\Lambda^{4}) \cap \ker d^{c}$
in $d\holda{k+1}(\Lambda^{4})$.
Hence $\harm^{5} \oplus d\holda{k+1}(\Lambda^{4}_{6})$ and
$\zfive_{k} \cap \ker d^{c}$ are both the \ltwoorth{} complement to
$W_{2}$ in $\zfive_{k}$, so they must be equal.
\end{proof}

\begin{defn}
Let $U \subseteq \ztre_{k} \times \ztwo_{k}$ be a small
neighbourhood of $(\Omega,\omega)$, and
\begin{gather*}
F : U \to W_{1} \times (\harm^{5} \oplus d\holda{k+1}(\Lambda^{4}_{6}))
\times \holda{k}(\Lambda^{6}), \;\: \\
(\beta, \gamma) \mapsto
(F_{1}(\beta), F_{2}(\beta, \gamma), F_{3}(\beta, \gamma)) ,
\end{gather*}
where $F_{1}(\beta)$ is defined by \textup{(\ref{fonedefeq})},
$F_{2}(\beta, \gamma) = P_{2}(\beta \wedge \gamma)$ with $P_{2}$ defined by
\textup{(\ref{p2defeq})}, and
$F_{3}(\beta, \gamma) = \quart \beta \wedge \hat\beta - \sixth \gamma^{3}$.
\end{defn}

\begin{prop}
Let $(\beta, \gamma)$ be a zero of $F$ sufficiently close to $(\Omega,\omega)$.
Then $(\beta,\gamma)$ is a \cystr.
\end{prop}

\begin{proof}
By proposition \ref{foneprop} $F_{1}(\beta) = 0$ implies that
$d\hat\beta = 0$.
Using lemma \ref{cy3charlem} it suffices to show that $\beta \wedge \gamma = 0$.
Since $\beta$ is a torsion-free \slstr{} it defines an integrable complex
structure~$J_{\beta}$. Let $d^{c}_{\beta}$ be the conjugate differential
with respect to $J_{\beta}$. Since $d\gamma = 0$ the 
$(3,0)+(0,3)$-part of $d^{c}_{\beta}\gamma$ with respect to $J_{\beta}$
vanishes. Thus
\begin{equation}
\label{dconeeq}
d^{c}_{\beta}(\beta \wedge \gamma) = - \beta \wedge d^{c}_{\beta}\gamma = 0 .
\end{equation}
Let $\eta = \beta \wedge \gamma$. 
$F_{2}(\beta,\gamma) = 0$ implies that
$\eta \in W_{2}$.
$d^{c} : \sob{1}(\Lambda^{5}) \to L^{2}(\Lambda^{6})$ is bounded below
transverse to its kernel. $W_{2} \cap \ker d^{c} = 0$ by proposition
\ref{dccharprop}, so there is a constant $A$ independent of
$\eta \in W_{2}$ such that
\[ \norm{d^{c}\eta}_{L^{2}} \geq A\norm{\eta}_{\sob{1}} . \]
The map $C^{1}(\Lambda^{3}) \times \sob{1}(\Lambda^{2})
\to L^{2}(\Lambda^{5}), \; (\beta, \eta) \mapsto d^{c}_{\beta}\eta$
is differentiable in $\beta$ and bounded linear in $\eta$, so there is a
constant $B$ such that for $\beta$ sufficiently close to $\Omega$
\[ \norm{(d^{c}-d^{c}_{\beta})\eta}_{L^{2}} \leq
B\norm{\beta-\Omega}_{C^{1}}\norm{\eta}_{\sob{1}} . \]
Hence
\begin{equation}
\label{dctwoeq}
\norm{d^{c}_{\beta}\eta}_{L^{2}} \geq \norm{d^{c}\eta}_{L^{2}}
- \norm{(d^{c}-d^{c}_{\beta})\eta}_{L^{2}} \geq
(A-B\norm{\beta-\Omega}_{C^{1}})\norm{\eta}_{\sob{1}} .
\end{equation}
Combining (\ref{dconeeq}) and (\ref{dctwoeq}) gives that if
$\norm{\beta-\Omega}_{C^{1}} < A/B$ then $\beta \wedge \gamma = 0$.
\end{proof}

\begin{prop}
$(DF)_{(\Omega,\omega)} : \ztre_{k} \times \ztwo_{k} \to
W_{1} \times (\harm^{5} \oplus d\holda{k+1}(\Lambda^{4}_{6})) \times \holda{k}(\Lambda^{6})$
is surjective.
\end{prop}

\begin{proof}
Suppose $(\chi_{1}, \chi_{2}, \chi_{3}) \in
W_{1} \times (\harm^{5} \oplus d\holda{k+1}(\Lambda^{4}_{6})) \times \holda{k}(\Lambda^{6})$.
By proposition \ref{foneprop} there is $\sigma \in \ztre_{k}$ such that
$(DF_{1})_{\Omega}\sigma = \chi_{1}$.
By proposition \ref{dccharprop} there is $\tau \in \ztwo_{k}$ such that
$\Omega \wedge \tau = \chi_{2} - P_{2}(\sigma \wedge \omega)$.
For $f \in \holda{k+2}(\bbr), \lambda \in \bbr$
\[ (DF_{3})_{(\Omega,\omega)}(0,dd^{c}f+\lambda\omega) =
-\half \omega^{2}\wedge (dd^{c}f + \lambda \omega) =
-\half (\triangle f + \lambda)\omega^{3} . \]
Since $\holda{k}(\Lambda^{6}) = 
\left( \bbr\oplus\triangle\holda{k+2}(\bbr) \right)\omega^{3}$
we can therefore find $f, \lambda$ such that
\[ (DF_{3})_{(\Omega,\omega)}(\sigma,\tau+dd^{c}f+\lambda\omega) = \chi_{3} . \]
Then $(DF)_{(\Omega,\omega)}(\sigma,\tau+dd^{c}f+\lambda\omega) =
(\chi_{1}, \chi_{2}, \chi_{3})$.
\end{proof}

Recall that the tangent space to the slice is
$\slt = (\harm^{3} \oplus W_{1}) \times \ztwo_{k}$, and that in definition
\ref{hsudef} we took $\harm_{SU}$ to be the harmonic tangents at
$(\Omega,\omega)$ to the space of $SU(3)$-structures.

\begin{prop}
The kernel of $(DF)_{(\Omega,\omega)}$ in $\slt$ is $\harm_{SU}$.
\end{prop}

\begin{proof}
It is obvious that $\harm_{SU}$ is contained in the kernel.
The projection from the kernel to the $\ztre_{k}$-component has image
contained in $\harm^{3} \subseteq \ztre_{k}$ and kernel contained in
$0 \times \harm^{2}_{8} \subseteq \slt$,
so by dimension-counting $\harm_{SU}$ is all of the kernel.
\end{proof}

As a consequence we have that $\slt$ really is
transverse to~$\orbt$, as claimed earlier.

\begin{cor}
\label{transcor}
$\slt \cap \orbt = 0$.
\end{cor}

\begin{proof}
If $\slt \cap \orbt$ were non-trivial the kernel of
$(DF)_{(\Omega,\omega)}$ in $\slt$ would contain some non-zero
exact forms.
\end{proof}

Now we can apply the implicit function theorem to $F$ to show that the
$\mathcal{Q}$ is a manifold with tangent space $\harm_{SU}$ at
$(\Omega,\omega)$.
Proposition \ref{cypremodprop} follows by elliptic regularity and
slice arguments.

\section{Hodge theory}
\label{hodgesec}

We wish to study the moduli space of torsion-free EAC \gstr s on a manifold
with cylindrical ends in terms of the projection (\ref{derhamprojeq})
to the de Rham cohomology.
In order to do this we require results about Hodge theory on EAC manifolds.
We also explain how the type decompositions of de Rham cohomology discussed in
subsection \ref{harmholsub} behave on EAC \gmfd s.

\subsection{Analysis of the Laplacian}
\label{lapsubsec}

We review some results that we shall need about analysis of elliptic
asymptotically translation-invariant operators on manifolds with
cylindrical ends, and explain how they can be applied to the Laplacian of an
EAC metric. For more detail about the analysis see e.g.
Lockhart \cite{lockhart87}, Lockhart and McOwen \cite{lockhart85},
and Maz'ya and Plamenevski\u{\i} \cite{mazya78}.

\begin{defn}
\label{wholddef}
Let $M^{n}$ be a manifold with cylindrical ends. If $g$ is an asymptotically
translation-invariant metric on $M$, $E$ is a vector bundle
associated to~$TM$, $\delta \in \bbr$ and $s$ is a section of $E$
define the \emph{Hölder norm with weight $\delta$}
(or $\holdad{k}$-norm) of $s$ in terms of the Hölder norm associated to $g$ by
\begin{equation*}
\norm{s}_{\holdad{k}(g)} = \norm{e^{\delta t}s}_{\holda{k}(g)} ,
\end{equation*}
where $t$ is the cylindrical coordinate on $M$.
Denote the space of sections of $E$ with finite $\holdad{k}$-norm
by $\holdad{k}(E)$.
\end{defn}

Up to equivalence, the weighted norms are independent of the choice of EAC
metric~$g$, and of the choice of $t$ on the compact piece $M_{0}$.
In particular, as topological vector spaces $\holdad{k}(E)$ are independent of
these choices.

We will want to use that $d$ and $d^{*}$ are formal adjoints in
integration by parts arguments. On a manifold with cylindrical ends
this is only justified if the rate of decay of the product is positive.

\begin{lem}
\label{ibplem}
Let $M^{n}$ be a manifold with cylindrical ends equipped with an asymptotically
translation-invariant metric.
Suppose that $\beta \in \holdax{k}{\delta_{1}}(\Lambda^{m}T^{*}M)$,
$\gamma \in \holdax{k}{\delta_{2}}(\Lambda^{m+1}T^{*}M)$ with $k \geq 1$ and
$\delta_{1} + \delta_{2} > 0$. Then
\[ \linner{d\beta,\gamma} \,=\, \linner{\beta,d^{*}\gamma} . \]
\end{lem}

Let $M$ be a manifold with cylindrical ends, $E, F$ vector bundles
associated to $TM$ and $A$ a linear smooth order $r$ differential operator
$\Gamma(E) \to \Gamma(F)$.
The restriction of $A$ to the cylindrical end $M_{\infty}$ can be written in
terms of the Levi-Civita connection of an arbitrary product metric
on $X \times \bbr$ as
\begin{equation}
\label{asycoeffeq}
A = \sum_{i = 0}^{r}a_{i}\nabla^{i}
\end{equation}
with coefficients $a_{i} \in C^{\infty}((TM)^{i} \otimes E^{*} \otimes F)$.
$A$ is said to be \emph{asymptotically translation-invariant} 
if the coefficients $a_{i}$~are. 
Then $A$ induces bounded linear maps
\begin{equation*}
\label{tiopeq}
A : \holdad{k+r}(E) \to \holdad{k}(F)
\end{equation*}
for any $\delta \in \bbr$.
One of the main results of \cite{lockhart85} is Theorem $6.2$, which states
that if $A$ is elliptic then these maps are Fredholm for all but a discrete
set of values of $\delta$,
and also relates the index for different values of $\delta$.
\cite[Theorem $7.4$]{lockhart85} is a corollary, which computes the index
of self-adjoint asymptotically translation-invariant elliptic operators for
small weights. This can be applied in particular to the Hodge Laplacian of
an asymptotically cylindrical metric, as in \cite[Section $3$]{lockhart87}.

\begin{prop}
\label{lapindprop}
Let $M$ be an asymptotically cylindrical manifold, and
$\epsilon_{1}$ the largest real such that
\begin{equation}
\label{lapmapeq}
\triangle : \holdax{k+2}{\pm\delta}(\Lambda^{m}T^{*}M) \to
\holdax{k}{\pm\delta}(\Lambda^{m}T^{*}M)
\end{equation}
is Fredholm for all $m$ and $0 < \delta < \epsilon_{1}$.
Then the index of \textup{(\ref{lapmapeq})} is $\mp (b^{m-1}(X) + b^{m}(X))$
for all $0 < \delta < \epsilon_{1}$.
\end{prop}

\begin{rmk}
Strictly speaking, the results in \cite{lockhart85} use weighted Sobolev
spaces rather than weighted Hölder spaces, but the arguments are the same
in both cases. See also \cite[Theorem $6.4$]{mazya78}.
\end{rmk}

\begin{lem}
\label{epsionectslem}
$\epsilon_{1}$ depends only on the asymptotic model $g_{X}+dt^{2}$
for the metric on $M$. Furthermore, $\epsilon_{1}$ is a lower semi-continuous
function of $g_{X}$ with respect to the $C^{1}$-norm.
\end{lem}

\begin{proof}
$\epsilon_{1}^{2}$ is in fact the smallest positive eigenvalue $\lambda_{1}$ of
the Hodge Laplacian $\triangle_{X}$ defined by $g_{X}$ on $\Omega^{*}(X)$.
To prove the proposition it therefore suffices to show that $\lambda_{1}$ is
lower semi-continuous in $g_{X}$.

Let $g, g'$ be smooth Riemannian metrics on $X$, $\triangle, \triangle'$
their Laplacians and $\lambda_{1}, \lambda_{1}'$ the smallest positive
eigenvalues of the Laplacians.
Let $T$ be the $L^{2}(g)$-orthogonal complement to $\ker \triangle$ in
$\holda{2}(\Lambda^{*}T^{*}X)$.
Then for any $\beta \in T$ with unit $L^{2}(g)$-norm
\[ \lambda_{1} \leq \; \lginner{\triangle \beta, \beta} \; =
\lgnorm{d\beta}^{2}+\lgnorm{d^{*}\beta}^{2} . \]
Since $d+d^{*}$ is an elliptic operator it gives a Fredholm map
$\sob{1}(\Lambda^{*}T^{*}X) \to L^{2}(\Lambda^{*}T^{*}X)$, so it is
bounded below transverse to its kernel. In other words, there is
a constant $C_{1}$ such that
$\norm{\beta}^{2}_{\sob{1}(g)} \leq C_{1}
\left(\lgnorm{d\beta}^{2}+\lgnorm{d^{*}\beta}^{2}\right)$
for any $\beta \in T \cap \sob{1}(\Lambda^{*}T^{*}X)$.

Let $e_{1}$ be an eigenvector of $\triangle'$ with eigenvalue $\lambda_{1}'$.
By Hodge theory for compact manifolds $\ker \triangle$ and $\ker \triangle'$
have the same dimension, so $(\ker \triangle' \oplus \bbr e_{1}) \cap T$ is
non-trivial. Hence
\begin{equation*}
\lambda_{1}' \geq
\frac{\lgpinner{\triangle' \beta, \beta}}
{\lgpinner{\beta,\beta}} =
\frac{\lgpnorm{d\beta}^{2} + \lgpnorm{d^{*'}\beta}^{2}}
{\lgpnorm{\beta}^{2}}
\end{equation*}
for some $\beta \in T$ with unit $L^{2}(g)$-norm. The RHS
depends differentiably on $g'$ (with respect to the $C^{1}(g)$-norm) and the
derivative at $g' = g$
can be estimated in terms of $\norm{\beta}^{2}_{\sob{1}(g)}$.
Therefore there is a constant $C_{2}$ (independent of $\beta$)
such that for any $g'$ close to $g$
\begin{multline*}
\lambda_{1}' \geq
\frac{\lgpnorm{d\beta}^{2} + \lgpnorm{d^{*'}\beta}^{2}}
{\lgpnorm{\beta}^{2}} \\
\geq \lgnorm{d\beta}^{2}+\lgnorm{d^{*}\beta}^{2}
- C_{2}\norm{g'-g}_{C^{1}(g)}\norm{\beta}^{2}_{\sob{1}(g)}
\geq \left(1 - \norm{g'-g}_{C^{1}(g)}C_{1}C_{2}\right)\lambda_{1} . \qedhere
\end{multline*}
\end{proof}

Now let $M$ be an EAC manifold with rate $\delta_{0}$ and cross-section $X$,
and assume that $0 < \delta < \min\{\epsilon_{1}, \delta_{0}\}$.
We fix some notation for various spaces of harmonic forms.

\begin{defn}
\label{harmdef}
Denote by
\begin{enumerate}
\item $\harm^{m}_{\pm}$ the space of harmonic $m$-forms in
$\holdax{k}{\pm\delta}(\Lambda^{m}T^{*}M)$,
\item $\harm^{m}_{0}$ the space of bounded harmonic $m$-forms on $M$,
\item $\harm^{m}_{\infty}$ the space of translation-invariant
harmonic $m$-forms on $X \times \bbr$,
\item $\harm^{m}_{X}$ the space of harmonic $m$-forms on $X$.
\end{enumerate}
\end{defn}

By elliptic regularity $\harm^{m}_{\pm}$ consists of smooth forms, and is
independent of $k$ for $k \geq 2$.
The computation of the index in proposition \ref{lapindprop}
involves proving that
\begin{equation}
\label{harminftyeq}
\harm^{m}_{\infty} = \{\psi + dt \wedge \tau : \psi \in \harm^{m}_{X},
\tau \in \harm^{m-1}_{X} \} .
\end{equation}
In particular the index of
\begin{equation}
\label{nonpadeq}
\triangle : \holdad{k+2}(\Lambda^{m}T^{*}M) \to \holdad{k}(\Lambda^{m}T^{*}M)
\end{equation}
is $-\dim \harm^{m}_{\infty}$ for small positive $\delta$.
Knowing the index of the Laplacian on weighted Hölder spaces allows us to
use an index-counting argument to deduce results about the kernel, and
a Hodge decomposition statement.

Let $i = b^{m}(X) + b^{m-1}(X)$. $\harm^{m}_{\infty}$ has dimension $i$,
and the index of $(\ref{nonpadeq})$ is $-i$.
$\triangle(\rho\harm^{m}_{\infty})$ and $\triangle(\rho t \harm^{m}_{\infty})$
consist of exponentially decaying forms.
Therefore
\begin{equation}
\label{padeq}
\triangle : \holdad{k+2}(\Lambda^{m}T^{*}M) \oplus \rho\harm^{m}_{\infty}
\oplus \rho t \harm^{m}_{\infty} \to \holdad{k}(\Lambda^{m}T^{*}M)
\end{equation}
is well-defined, and its index is $+i$. (\ref{nonpadeq}) has kernel
$\harm^{m}_{+}$, and by integration by parts the image contained in the
orthogonal complement of $\harm^{m}_{-}$. (\ref{padeq}) has kernel contained
in $\harm^{m}_{-}$ and image contained in the orthogonal complement of
$\harm^{m}_{+}$. Hence
{\setlength\arraycolsep{2pt}
\begin{eqnarray*}
i & \leq & \dim \harm^{m}_{-} - \dim \harm^{m}_{+} , \\
-i & \leq & \dim \harm^{m}_{+} - \dim \harm^{m}_{-} .
\end{eqnarray*}}
Since equality holds the kernel of (\ref{padeq}) is exactly $\harm^{m}_{-}$.

\begin{prop}
\label{explodeprop}
Let $M$ be an exponentially asymptotically cylindrical manifold with rate
$\delta_{0}$, $k \geq 0$ and $0 < \delta < \min\{\epsilon_{1}, \delta_{0}\}$.
Then
\begin{gather*}
\harm^{m}_{-} \subseteq \holdad{k}(\Lambda^{m}T^{*}M) \oplus
\rho\harm^{m}_{\infty} \oplus \rho t \harm^{m}_{\infty} , \\
\harm^{m}_{0} \subseteq \holdad{k}(\Lambda^{m}T^{*}M) \oplus
\rho\harm^{m}_{\infty} ,
\end{gather*}
and $\harm^{m}_{+}$ is precisely the space of $L^{2}$-integrable harmonic forms
on $M$.
\end{prop}

Also the image of (\ref{padeq}) is exactly the \ltwoorth{} complement of
$\harm^{m}_{+}$ in $\holdad{k}(\Lambda^{m}T^{*}M)$.
We can use this to prove an EAC analogue of the Hodge decomposition
for compact manifolds:
\begin{equation}
\label{hodgedecompeq}
\holdad{k}(\Lambda^{m}T^{*}M) = \harm^{m}_{+} \oplus
\holdad{k}[d\Lambda^{m-1}T^{*}M] \oplus \holdad{k}[d^{*}\Lambda^{m+1}T^{*}M] ,
\end{equation}
where we let $\holdad{k}[d\Lambda^{m-1}T^{*}M]$ and
$\holdad{k}[d^{*}\Lambda^{m+1}T^{*}M]$
denote the subspaces of $\holdad{k}(\Lambda^{m}T^{*}M)$ consisting of exact and
coexact forms, respectively.

\subsection{Hodge theory on EAC manifolds}
\label{hodgesub}

In this subsection we outline the correspondence between bounded harmonic forms
and the de Rham cohomology on an oriented EAC Riemannian manifold $M^{n}$ with
cross-section~$X^{n-1}$. For `exact $b$-manifolds', which can be considered
as a subclass of EAC manifolds, the type of results we need can be found in
Section $6.4$ of Melrose \cite{melrose94}.
The arguments there carry over unchanged to the EAC case.

\begin{defn}
Let $\harm^{m}_{E}, \harm^{m}_{E^{*}} \subseteq \harm^{m}_{0}$
denote the spaces of bounded exact and coexact harmonic forms respectively
on $M$, and
\[ \harm^{m}_{abs} = \harm^{m}_{+} \oplus \harm^{m}_{E^{*}}, \;\;
\harm^{m}_{rel} = \harm^{m}_{+} \oplus \harm^{m}_{E} . \]
\end{defn}

\noindent The next theorem is part of \cite[Theorem $6.18$]{melrose94}.

\begin{thm}
\label{absderhamthm}
Let $M$ be an oriented EAC manifold. The natural projection map
$\pi_{H} : \harm^{m}_{abs} \to H^{m}(M)$ is an isomorphism.
\end{thm}

Since the kernel of $\pi_{H} : \harm^{m}_{0} \to H^{3}(M)$ is $\harm^{m}_{E}$
the theorem implies that
\begin{equation}
\harm^{m}_{0} = \harm^{m}_{+} \oplus \harm^{m}_{E^{*}} \oplus \harm^{m}_{E} .
\end{equation}
By proposition \ref{explodeprop} any $\alpha \in \harm^{m}_{0}$ is
exponentially asymptotic to some $B({\alpha}) \in \harm^{m}_{\infty}$.
By (\ref{harminftyeq})
$\harm^{m}_{\infty} \cong \harm^{m}_{X} \oplus dt \wedge \harm^{m-1}_{X}$,
so we get a boundary map
\[ B : \harm^{m}_{0} \to \harm^{m}_{X} \oplus dt \wedge \harm^{m-1}_{X}, \;\:
\alpha \mapsto B_{a}(\alpha) + dt \wedge B_{e}(\alpha) . \]
It is easy to see that for $\alpha \in \harm^{m}_{0}$ the pull-back map $j^{*}$
in the long exact sequence for relative cohomology (\ref{relexacteq}) acts as
\begin{equation}
\label{jbaeq}
j^{*}([\alpha]) = [B_{a}(\alpha)] \in H^{m}(X) ,
\end{equation}
and it follows that $\harm^{m}_{rel} \subseteq \ker B_{a}$. Applying the Hodge
star shows that also $\harm^{m}_{abs} \subseteq \ker B_{e}$.
Therefore $B_{a}$ is injective on
$\harm^{m}_{E^{*}}$ and $0$ on $\harm^{m}_{rel}$, while
$B_{e}$ is injective on $\harm^{m}_{E}$ and $0$ on~$\harm^{m}_{abs}$.

As a corollary of theorem \ref{absderhamthm} we can determine that
the image of the space $\harm^{m}_{+}$ of $L^{2}$ harmonic forms in the
de Rham cohomology $H^{m}(M)$ is precisely the subspace $H^{m}_{0}(M)$ of
compactly supported classes. This result appears as 
e.g. \cite[Proposition $4.9$]{atiyah75}, 
\cite[Theorems $7.6$ and $7.9$]{lockhart87},
and \cite[Proposition $6.14$]{melrose94}.

\begin{thm}
\label{derhamthm}
Let $M$ be an oriented EAC manifold. Then
$\pi_{H} : \harm^{m}_{+} \to H^{m}(M)$ is an isomorphism onto
$H^{m}_{0}(M)$.
\end{thm}

\begin{proof}
$\harm^{m}_{+}$ is $\ker B_{a}$ in $\harm^{m}_{abs}$, and it follows from
theorem \ref{absderhamthm} that it is mapped isomorphically to
$H^{m}_{0}(M) = \ker j^{*} \subseteq H^{m}(M)$.
\end{proof}

\begin{defn}
Let $\calm{A}{m} = B_{a}(\harm^{m}_{0}) \subseteq \harm^{m}_{X}, \, \:
\calm{E}{m} = B_{e}(\harm^{m+1}_{0}) \subseteq \harm^{m}_{X}$, and
let $A^{m}, E^{m}$ be the subspaces of $H^{m}(X)$ that they represent.
\end{defn}
$A^{m}$ is of course just the image $j^{*}(H^{m}(M)) \subseteq H^{m}(X)$.
The Hodge star on $M$ identifies $\harm^{m}_{abs}$
and $\harm^{m-n}_{rel}$. If $\beta \in \harm^{m}_{0}$ then
$B_{e}(*\beta) = *B_{a}(\beta)$. Therefore the Hodge star on $X$ identifies
$\mathcal{A}^{m}$ with $\mathcal{E}^{n-m-1}$, and
$A^{m}$ with $E^{n-m-1}$. \cite[Lemma $6.15$]{melrose94} implies

\begin{prop}
\label{aeprop}
Let $M^{n}$ be an oriented EAC manifold. Then
\[ \harm^{m}_{X} = \mathcal{A}^{m} \oplus \mathcal{E}^{m} \]
is an orthogonal direct sum.
\end{prop}

Finally, we observe

\begin{cor}
\label{degonehodgecor}
Let $M^{n}$ be an oriented EAC manifold which has a single end
(i.e.~the cross-section $X$ is connected). Then
$\inc : H^{1}_{cpt}(M) \to H^{1}(M)$ is injective.
$\harm^{1}_{E} = 0$, and $\harm^{1}_{0} \to H^{1}(M)$ is
an isomorphism.
\end{cor}

\begin{proof}
Consider the start of the long exact sequence for relative cohomology
\[ H^{0}_{cpt}(M) \to H^{0}(M) \to H^{0}(X) \stackrel{\partial}{\to}
H^{1}_{cpt}(M) \stackrel{\inc}{\to} H^{1}(M) . \]
The dimensions of the first three terms are $0, 1$, and $1$, so $\partial = 0$,
and thus $\inc$ is injective.
\[ \harm^{1}_{E} \cong \mathcal{E}^{0} \cong E^{0} \cong \partial(H^{0}(X)) =
\ker e \subseteq H^{1}_{cpt}(M) , \]
so the result follows.
\end{proof}

\subsection{Hodge theory of EAC $G_{2}$-manifolds}
\label{hodgeg2sub}

Let $M^{7}$ be an EAC \gmfd, with \gstr{} $\varphi$ asymptotic
to $\Omega + dt \wedge \omega$.
$(\Omega,\omega)$ is a Calabi-Yau structure on the cross-section $X$.
Maps such as
$\Omega^{2}(X) \to \Omega^{4}(X), \sigma \mapsto \sigma \wedge \omega$
are $SU(3)$-equivariant, so by proposition \ref{hodgelapprop} induce
maps between type components of the spaces of harmonic forms.
In this subsection we consider the relation between the type decompositions
and the decomposition in proposition \ref{aeprop}.

By remark \ref{cyldualrmk} $*\varphi$ is asymptotic to
$\half \omega^{2} - dt \wedge \hat\Omega$, where $\hat\Omega$ is the unique
$3$-form on $X$ such that $\Omega + i \hat\Omega$ has type $(3,0)$ as discussed
in subsection \ref{cydefsub}.

\begin{lem}
\label{aemapslem}
If $\tau \in \mathcal{E}^{m}$ then
$\tau \wedge \Omega \in \mathcal{E}^{m+3}$ and
$\tau \wedge \half \omega^{2} \in \mathcal{E}^{m+4}$.

If $\sigma \in \mathcal{A}^{m}$ then
$\sigma \wedge \Omega \in \mathcal{A}^{m+3}$,
$\sigma \wedge \half \omega^{2} \in \mathcal{A}^{m+4}$,
$\sigma \wedge \omega \in \mathcal{E}^{m+2}$ and
$\sigma \wedge \hat \Omega \in \mathcal{E}^{m+3}$.
\end{lem}

\begin{proof}
If $\chi \in \harm^{m+1}_{E}$ with $B_{e}\chi = \tau$ then
\begin{gather*}
\chi \wedge \varphi \in \harm^{m+4}_{0} \Rightarrow
dt \wedge \tau \wedge \Omega = B(\chi \wedge \varphi)
\in dt \wedge \mathcal{E}^{m+3} , \\
\chi \wedge *\varphi \in \harm^{m+5}_{0} \Rightarrow
dt \wedge \tau \wedge \half \omega^{2} = B(\chi \wedge *\varphi) \in
dt \wedge \mathcal{E}^{m+4} .
\end{gather*}
If $\chi \in \harm^{m}_{0}$ with $B_{a}\chi = \sigma$ then
\begin{gather*}
\chi \wedge \varphi \in \harm^{m+3}_{0} \Rightarrow
\sigma \wedge \Omega + dt \wedge \sigma \wedge \omega = B(\chi \wedge \varphi)
\in \mathcal{A}^{m+3} \oplus dt \wedge \mathcal{E}^{m+2} , \\
\chi \wedge *\varphi \in \harm^{m+4}_{0} \Rightarrow
\sigma \wedge \half \omega^{2} + dt \wedge \sigma \wedge \hat \Omega =
B(\chi \wedge *\varphi) \in
\mathcal{A}^{m+4} \oplus dt \wedge \mathcal{E}^{m+3} . \qedhere
\end{gather*}
\end{proof}

Hodge theory for compact manifolds allows us to identify the de Rham cohomology
of $X$ with the harmonic $m$-forms on $X$. The $L^{2}$-inner product on
$\harm^{m}_{X}$ therefore gives an inner product on $H^{m}(X)$, and
the Hodge star $* : \harm^{m}_{X} \to \harm^{6-m}_{X}$ gives isomorphisms
\[ * : H^{m}(X) \to H^{6-m}(X) . \]
This map is the composition of the metric isomorphism
$H^{m}(X) \cong (H^{m}(X))^{*}$ with Poincar\'e duality
$(H^{m}(X))^{*} \cong H^{6-m}(X)$.
Proposition \ref{aeprop} implies that there is an orthogonal direct sum
\begin{equation*}
\label{aeeq}
H^{m}(X) = A^{m} \oplus E^{m} ,
\end{equation*}
where $A^{m} = j^{*}(H^{m}(M))$ and $E^{m} = *A^{6-m}$.
Let $A^{2}_{6} = A^{2} \cap H^{2}_{6}(X)$ etc.

\begin{prop}
\label{aeisoprop}
Let $M$ be an EAC \gmfd{} with cross-section $X$. Then
\[ H^{2}_{6}(X) = A^{2}_{6} \oplus E^{2}_{6}, \;\;
H^{4}_{6}(X) = A^{4}_{6} \oplus E^{4}_{6} , \]
and the sums are orthogonal. Furthermore
\begin{enumerate}
\item \label{hodgeitem}
$H^{2}_{6}(X) \to H^{4}_{6}(X), [\alpha] \mapsto *[\alpha]$ maps
$A^{2}_{6} \to E^{4}_{6}, E^{2}_{6} \to A^{4}_{6}$,

\item \label{omitem}
$H^{1}(X) \to H^{4}_{6}(X), [\alpha] \mapsto
[\alpha] \hcup [\Omega]$ maps $A^{1} \to A^{4}_{6}, E^{1} \to E^{4}_{6}$,

\item \label{om2item}
$H^{1}(X) \to H^{5}(X), [\alpha] \mapsto
[\alpha] \hcup [\half \omega^{2}]$
maps $A^{1} \to A^{5}, E^{1} \to E^{5}$.
\end{enumerate}
\end{prop}

\begin{proof}
\ref{hodgeitem} is obvious, since $*$ maps $A^{m} \leftrightarrow E^{6-m}$.

$[\alpha] \mapsto [\alpha] \hcup [\Omega]$ is a bijection
$H^{1}(X) \to H^{4}_{6}(X)$, and it maps $A^{1} \hookrightarrow A^{4}$,
$E^{1} \hookrightarrow E^{4}$ by lemma \ref{aemapslem}. Thus \ref{omitem}.
It follows that that $A^{1} \to A^{4}_{6}, E^{1} \to E^{4}_{6}$
are both surjective, and that
$H^{4}_{6}(X)$ splits as $A^{4}_{6} \oplus E^{4}_{6}$.
$H^{2}_{6}(X)$ splits too by~\ref{hodgeitem}.

\ref{om2item} follows from lemma \ref{aemapslem} in the same way.
\end{proof}

When $X$ is Kähler the complex structure $J$ is parallel,
and the Hodge Laplacian on forms commutes with the action of $J$.
Thus $\harm^{1}_{X}$ inherits a complex structure from $X$ in the Kähler case,
and the inner product on $\harm^{1}_{X}$ is Hermitian.
We identify $H^{1}(X) \leftrightarrow \harm^{1}_{X}$ as usual to make
$H^{1}(X)$ a Hermitian vector space.

\begin{prop}
\label{realsubprop}
Let $M$ be an EAC \gmfd{} with cross-section $X$. Then the pull-back
$j^{*} : H^{1}(M) \hookrightarrow H^{1}(X)$ embeds $H^{1}(M)$ as a
Lagrangian subspace of $H^{1}(X)$ with the symplectic form
$<\cdot, J\cdot>$.
In particular $b^{1}(M) = \half b^{1}(X)$.
\end{prop}

\begin{proof}
By proposition \ref{cylbochnerprop} below $j^{*}$ maps $H^{1}(M)$
isomorphically to its image $A^{1}$.
The complex structure on $H^{1}(X)$ can be expressed as
\[ J[\alpha] = *([\alpha] \hcup [\half \omega^{2}]) . \]
Thus $J$ maps $A^{1}$ to its orthogonal complement $E^{1}$ by
proposition \ref{aeisoprop}.
\end{proof}

On compact Ricci-flat manifolds harmonic $1$-forms are parallel, and this
can be generalised to the EAC case.

\begin{prop}
\label{cylbochnerprop}
If $M$ is a Ricci-flat EAC manifold then $\harm^{1}_{0}$ is the space of
parallel $1$-forms on $M$. In particular $\harm^{1}_{+} = 0$, and
$j^{*} : H^{1}(M) \to H^{1}(X)$ is injective.
\end{prop}

\begin{proof}
This is proved by a standard `Bochner argument'.
For a 1-form $\phi$ the Weitzenböck formula
\cite[(1.155)]{besse87} gives
\begin{equation}
\label{weitzeneq}
\triangle \phi = \nabla^{*}\nabla \phi .
\end{equation}
It follows that any parallel $1$-form $\phi$ is harmonic, and
parallel forms are of course bounded.
To show that any bounded harmonic form is parallel we use (\ref{weitzeneq})
together with an integration by parts argument.

$\harm^{1}_{+}$ consists of parallel decaying forms, so is $0$. By theorem
\ref{derhamthm} the kernel $H^{1}_{0}(M)$ of $j^{*} : H^{1}(M) \to H^{1}(X)$ is
represented by $\harm^{1}_{+}$.
\end{proof}

\section{Deformation theory of EAC $G_{2}$-manifolds}
\label{maincylsec}

In this section we construct the moduli space of exponentially
asymptotically cylindrical torsion-free \gstr s on a manifold with a
cylindrical end, and study its local properties.

Throughout the section $M^{7}$ is a connected oriented manifold with
cylindrical ends, $X^{6}$ is its cross-section, $t$ the cylindrical
coordinate on $M$ and $\rho$ a cut-off function for the cylinder on $M$.
We abbreviate $\Lambda^{m}T^{*}M$ to~$\Lambda^{m}$.

\subsection{Proof outline}
\label{outlinesub}

We now set out to prove the main theorem \ref{maincylthm}. The argument is
a generalisation of that used by Hitchin in \cite{hitchin00} to construct
the moduli space of torsion-free \gstr s on a compact manifold. Hitchin shows
that a \gstr{} is torsion-free if and only if it is a critical point of the
volume functional $\varphi \mapsto Vol(\varphi)$ ($Vol(\varphi)$ is the total
volume of the metric defined by $\varphi$). An appropriate modification of
the map $\varphi \mapsto D_{\varphi}Vol$ has surjective derivative, and the
implicit function theorem can be applied to find pre-moduli spaces of
torsion-free \gstr s which are manifolds.

In the EAC case it is not as natural to consider the volume functional,
but we can nevertheless adapt the steps in Hitchin's proof to find
pre-moduli spaces which are manifolds, and then apply slice arguments
to show that the moduli space is a manifold. Like in e.g. Kovalev
\cite{kovalev05}
(which studies deformations of EAC Calabi-Yau metrics) we can simplify
the slice by understanding the deformations of the boundary -- in this
case we use the deformation theory of compact Calabi-Yau $3$-folds described
in section $5$.

As an intermediate step in the proof of theorem \ref{maincylthm} we construct
moduli spaces of torsion-free exponentially cylindrical \gstr s with some
fixed rate $\delta > 0$. As before we let $\cald{X}$ be the space of
torsion-free EAC \gstr s with rate $\delta$.
Each $\varphi \in \cald{X}$ defines an EAC metric, and hence a parameter
$\epsilon_{1}(\varphi)$ (cf.~proposition \ref{lapindprop}). When we study
a neighbourhood of $\varphi$ we need to assume that
$\delta < \epsilon_{1}(\varphi)$ in order to apply analysis results.
We therefore let
\begin{equation*}
\calxa = \{ \varphi \in \cald{X} : \epsilon_{1}(\varphi) > \delta \} .
\end{equation*}
$\epsilon_{1}$ depends lower semi-continuously on the
asymptotic model by lemma~\ref{epsionectslem}, so
$\calxa$ is an open subset of $\cald{X}$.
Let $\cald{D}$ be the group of EAC 
diffeomorphisms of $M$ with rate~$\delta$. The rate $\delta$ moduli space
that we study is \mbox{$\cald{M} = \calxa/\cald{D}$}.

Given $\varphi \in \calxa$ let 
$\Omega + dt \wedge \omega = B(\varphi)$, i.e. the asymptotic limit of
$\varphi$.
If we identify $\Omega + dt \wedge \omega$ with the pair $(\Omega,\omega)$ then
by proposition \ref{cylg2prop}
$\Omega + dt \wedge \omega$ defines a \cystr{} on $X$.
By proposition \ref{cypremodprop} there is a pre-moduli space $\mathcal{Q}$
of Calabi-Yau structures near $\Omega + dt \wedge \omega$. Let
\begin{equation*}
\calq{X} = \{ \psi \in \calxa : B(\psi) \in \mathcal{Q} \} ,
\end{equation*}
and let $\calq{D} \subseteq \cald{D}$ be the subgroup of diffeomorphisms
asymptotic to automorphisms of the cylindrical \gstr{}
$\Omega + dt \wedge \omega$. By proposition \ref{autinvprop} $\calq{D}$ acts
on $\calq{X}$, and we will see that $\calq{X}/\calq{D}$ maps homeomorphically
to an open subset of $\cald{M}$.

We use slice arguments to study a neighbourhood of $\varphi\calq{D}$ in
$\calq{X}/\calq{D}$. In order to be able to apply analysis results we need
to use Banach spaces of forms, so we work with weighted Hölder $\holdad{k}$
spaces, for some fixed $k \geq 1, \alpha \in (0,1)$.
Note that the boundary values of elements of $\calq{X}$ must lie in
\begin{equation}
\label{qaeq}
\mathcal{Q}_{A} = \{\Omega' + dt \wedge \omega' \in \mathcal{Q}
: [\Omega'] \in A^{3}, [\omega'^{2}] \in A^{4} \} ,
\end{equation}
where $A^{m}$ is the image of $j^{*} : H^{m}(M) \to H^{m}(X)$ (cf. discussion
before theorem \ref{modsubmersethm}). We use the cut-off function $\rho$ to
consider $\rho\mathcal{Q}_{A}$ as a subspace of smooth asymptotically
translation-invariant $3$-forms supported on the cylinder of $M$, and let
\[ \ztreq \subseteq \holdad{k}(\Lambda^{3}) + \rho \mathcal{Q}_{A} \]
be the subspace of closed forms. Then $\calq{X} \hookrightarrow \ztreq$
continuously.
The main technical steps in the construction of the pre-moduli space near
$\varphi$ are to
\begin{enumerate}
\item show that $\mathcal{Q}_{A}$ is a submanifold of $\mathcal{Q}$
(proposition \ref{qasubmfdprop}), so that $\ztreq$ is a manifold,
\item find a complement $\slt$ in $T_{\varphi}\ztreq$ for the tangent
space to the $\calq{D}$-orbit at $\varphi$
(proposition \ref{sltprop}),
and pick a submanifold $\mathcal{S} \subseteq \ztreq$ with
$T_{\varphi}\mathcal{S} = \slt$,
\item
show that the space of torsion-free
\gstr s $\cald{R} \subseteq \mathcal{S}$ is a submanifold
(proposition \ref{premodprop}),
\item show that the elements of $\cald{R}$ are smooth and exponentially
asymptotically translation-invariant with rate $\delta$
(proposition \ref{cylregprop}).
\end{enumerate}

In subsection \ref{modconstrsub} we then provide the slice arguments that
show that $\cald{R}$ is homeomorphic to a neighbourhood in $\calq{X}/\calq{D}$,
and therefore in~$\cald{M}$. Hence $\cald{M}$ is a manifold for any
$\delta > 0$. We then show that $\cald{M}$ is homeomorphic to an open
subset of $\calp{M}$ for any $\delta > 0$, and deduce that $\calp{M}$ is
a manifold.

\begin{rmk}
The last step means that if $\varphi \in \calp{X}$ is EAC with rate
$\delta_{0}(\varphi)$ then for any
$0 < \delta < \min \{\delta_{0}(\varphi), \epsilon_{1}(\varphi)\}$
the pre-moduli
space $\cald{R}$ gives a chart near $\varphi$ not only in $\cald{M}$, but also
in $\mathcal{M}_{\delta'}$ for any $\delta' > \delta$, and in $\calp{M}$.
In other words, $\cald{R}$ is essentially independent of $\delta$ if $\delta$
is chosen sufficiently small.
\end{rmk}

\subsection{The boundary values}

As explained above we are restricting our attention to determining the
space of torsion-free \gstr s in $\ztreq$, whose boundary values
lie in a space $\mathcal{Q}_{A}$. In order to make sense of this we first of
all need to know that $\mathcal{Q}_{A}$ is a manifold. We show that here, and
in the process essentially prove theorem \ref{modsubmersethm}.

Let $X^{6}$ be the cross-section of an EAC \gmfd{} $M^{7}$, and
$(\Omega, \omega)$ a Calabi-Yau structure on $X$ defined by the limit
of a torsion-free EAC \gstr{} on $M$.
Let $\mathcal{Q}$ be the pre-moduli space of Calabi-Yau structures near
$(\Omega, \omega)$, and equivalently to (\ref{qaeq}) define
\begin{equation*}
\mathcal{Q}_{A} = \{(\Omega', \omega') \in \mathcal{Q}
: [\Omega'] \in A^{3}, [\omega'^{2}] \in A^{4} \} .
\end{equation*}
Since $\mathcal{Q}$ is diffeomorphic to a neighbourhood in the moduli space
$\mathcal{N}$ of \cystr s on~$X$, $\mathcal{Q}_{A}$ is homeomorphic to a
neighbourhood in the subspace $\mathcal{N}_{A}$ of classes which stand a chance
of being the boundary values of EAC torsion-free \gstr s, as discussed before
theorem \ref{modsubmersethm}.

Recall that by proposition \ref{cypremodprop} the tangent space at
$(\Omega, \omega)$ to the pre-moduli space $\mathcal{Q}$ is the space
$\harm_{SU}$ of harmonic tangents to the $SU(3)$-structures.
As before let $E^{m} \subseteq H^{m}(X)$ be the orthogonal complement of
$A^{m}$, and let $\mathcal{A}^{m}, \mathcal{E}^{m} \subseteq \harm^{m}_{X}$
denote the respective spaces of harmonic representatives.
By lemma \ref{aemapslem} $\tau \mapsto \omega \wedge \tau$ maps
$\mathcal{A}^{2} \to \mathcal{E}^{4}$ and $\mathcal{E}^{2} \to \mathcal{A}^{4}$.
Hence the linearisation of the condition $[\omega'^{2}] \in A^{4}$ is
$[\tau] \in E^{2}$, and we would expect the tangent space to $\mathcal{Q}_{A}$
at $(\Omega,\omega)$ to be
\begin{equation*}
\harm_{SU,A} =
\{ (\sigma, \tau) \in \harm_{SU}
: \sigma \in \mathcal{A}^{3}, \tau \in \mathcal{E}^{2} \} .
\end{equation*}

\begin{prop}
\label{qasubmfdprop}
Let $(\Omega,\omega)$ be the Calabi-Yau structure induced on the cross-section
$X^{6}$ of an EAC \gmfd{} $M^{7}$, and $\mathcal{Q}$ the
pre-moduli space of \cystr s near $(\Omega,\omega)$.
Then $\mathcal{Q}_{A} \subseteq \mathcal{Q}$ is a submanifold, and
\[ T_{(\Omega, \omega)}\mathcal{Q}_{A} = \harm_{SU,A} . \]
\end{prop}

\begin{proof}
The map $\mathcal{Q} \to H^{3}(X)$ is a submersion, so
\[ \mathcal{Q}' =
\{(\Omega',\omega') \in \mathcal{Q} : [\Omega'] \in A^{3} \} \]
is a submanifold of $\mathcal{Q}$.
By proposition \ref{aeisoprop}
\[ H^{4}(X) = A^{4}_{1} \oplus A^{4}_{6} \oplus A^{4}_{8}
\oplus E^{4}_{6} \oplus E^{4}_{8} , \]
where $A^{4}_{6} = H^{4}_{6}(X) \cap A^{4}$ etc.
Let
\[ P_{E8} : H^{4}(X) \to E^{4}_{8} \]
be the orthogonal projection. For
$(\Omega', \omega') \in \mathcal{Q}'$ let ${E^{m}}'$ be the orthogonal
complement of $A^{m}$ in $H^{m}(X)$ with respect to the metric defined by
$(\Omega', \omega')$, and let
$P_{A'} : H^{m}(X) \to A^{m}$ and $P_{E'} : H^{m}(X) \to {E^{m}}'$ be the
projections.
Let
\[ F : \mathcal{Q'} \to E^{4}_{8}, \;\:
(\Omega',\omega') \mapsto P_{E8}P_{E'}[\omega'\wedge\omega'] . \]
We prove that $\mathcal{Q}_{A}$ is a submanifold of $\mathcal{Q}'$ by
showing that it is the zero set of~$F$, and that $F$ has surjective derivative
at $(\Omega,\omega)$.

Suppose $F(\Omega',\omega') = 0$, and let $a = P_{E'}[\omega'\wedge\omega']$.
Write $a = b + c$, with $b \in A^{4}$,
$c \in E^{4}$. $P_{E8}a = 0 \Rightarrow \pi_{8}c = 0$,
so $c \in E^{4}_{6}$.
Since $E^{1} \to E^{4}_{6}, \; v \mapsto [\Omega] \hcup v$ is an isomorphism
$c = [\Omega] \hcup v$ for some $v \in E^{1}$. In the inner product
$\inner{\;,\;}'$ on $H^{*}(X)$ defined by $(\Omega',\omega')$
\begin{multline*}
\inner{a,a}' \;=\; \inner{a,[\Omega] \hcup v}' \;=\;
\inner{a, [\Omega] \hcup v - [\Omega'] \hcup P_{E'}v}' \leq
\norm{a}'(\norm{[\Omega-\Omega']}'\norm{v}' + \norm{P_{A'}v}') .
\end{multline*}
The RHS can be estimated by $\norm{[\Omega-\Omega']}(\norm{a}')^{2}$ for
$(\Omega',\omega')$ close to $(\Omega,\omega)$.
Hence for $(\Omega',\omega')$ sufficiently close to $(\Omega,\omega)$
\[ F(\Omega',\omega') = 0 \Rightarrow P_{E'}[\omega'^{2}] = 0
\Rightarrow [\omega'^{2}] \in A^{4} . \]
So $\mathcal{Q}_{A} \subseteq \mathcal{Q}'$ is the zero set of $F$.
It remains to show that $F$ has surjective derivative.
If $(\sigma, \tau) \in (\mathcal{A}^{3} \times \harm^{2}_{X}) \cap \harm_{SU} =
T_{(\Omega,\omega)}\mathcal{Q}'$ then since $[\omega^{2}] \in A^{4}$
\[DF_{(\Omega,\omega)}(\sigma, \tau) = P_{E8}P_{E}(2[\omega \wedge \tau])
= 2P_{E8}[\omega \wedge \tau] . \]
Since $\mathcal{A}^{2}_{8} \to \mathcal{E}^{4}_{8}, \;
\tau \mapsto \omega \wedge \tau$ is a bijection
the derivative maps $0 \times \mathcal{A}^{2}_{8} \subseteq
T_{(\Omega,\omega)}\mathcal{Q}'$ onto $E^{4}_{8}$.

By the implicit function theorem $\mathcal{Q}_{A}$ is a manifold,
and the tangent space at
$(\Omega,\omega)$ is
\[ \ker DF_{(\Omega,\omega)} = \harm_{SU,A} . \qedhere \]
\end{proof}

\begin{cor}
\label{hsuacor}
$\harm_{SU,A} \to \mathcal{A}^{3}, \; (\sigma,\tau) \mapsto \sigma$
is surjective, with kernel $0 \times \mathcal{E}^{2}_{8}$.
\end{cor}

\begin{proof}
The last part of the proof of the proposition actually shows that 
$\mathcal{Q}_{A} \to A^{3}$ is a submersion, so
$\harm_{SU,A} \to \mathcal{A}^{3}$ is surjective. This could also be deduced
from lemma \ref{aemapslem}. By definition of $\harm_{SU}$ the kernel
consists of $(0, \tau) \in 0 \times \mathcal{E}^{2}$ satisfying the conditions
(\ref{leqs}), which reduce to $\pi_{1}\tau = \pi_{6}\tau = 0$.
\end{proof}

Proposition \ref{qasubmfdprop} implies directly that a neighbourhood
of the image of $B : \calp{M} \to \mathcal{N}_{A}$ is a manifold.
The rest of theorem \ref{modsubmersethm} follows too, once we have proved the
main result that $\calp{M}$ is a manifold.
We will return to this in subsection \ref{propertiessub}.

\subsection{The Dirac operator}

We will use Fredholm properties of the Dirac operator associated to a \gstr{}
$\varphi$ to obtain a direct sum decomposition (proposition \ref{sltprop}).
The required properties
of the Dirac operator can conveniently be deduced from the properties
of the Laplacian discussed in subsection \ref{lapsubsec}.

Since $G_{2} \subseteq SO(7)$ is simply connected it can be regarded as a
subgroup of $Spin(7)$, so a \gstr{} on a manifold $M^{7}$ induces a
spin structure (in fact a converse also holds: an oriented manifold $M^{7}$
admits \gstr s if and only if it admits a spin structure,
cf.~\cite[Remark $3$]{bryant03}). The spin structure defines
a spinor bundle $S$, and the Dirac operator
\[ \dirac : \Gamma(S) \to \Gamma(S) . \]
Recall that if $\bbr^{7}$ is identified with the imaginary part of the
octonions $\bbo$ then $G_{2}$ corresponds to the group of normed-algebra
automorphisms of~$\bbo$. Indeed, the $3$-form
$\varphi_{0} \in \Lambda^{3}(\bbr^{7})^{*}$ stabilised
by $G_{2}$ can be written in terms of the octonion multiplication
$\cdot$ and the inner product $<,\!>$ as
\[ \varphi_{0}(a, b, c) =\; \inner{a \cdot b, c} . \]
The octonion multiplication by $\bbr^{7}$ on $\bbo$ induces a representation
of the Clifford algebra $Cl(\bbr^{7})$ on $\bbo$. Hence
$Spin(7) \subset Cl(\bbr^{7})$ also acts on $\bbo$, and this identifies
$\bbo$ with the spin representation of $Spin(7)$. This identification is
in fact $G_{2}$-equivariant (cf. \cite[page 122]{harvey82}).
In particular, the restriction of the spin representation to
$G_{2}$ is isomorphic to
$\bbr^{1} \oplus \bbr^{7}$, the direct sum of the trivial and standard
representations of~$G_{2}$.

Hence on a manifold $M^{7}$ with a \gstr{} $\varphi$ the spinor bundle
$S$ is isomorphic to
$\Lambda^{0}T^{*}M \oplus \Lambda^{1}T^{*}M$.
Under this identification Clifford multiplication by some
$\alpha \in \Omega^{1}(M)$ (which comes from the
octonion multiplication) corresponds to
\begin{equation*}
\Omega^{0}(M) \oplus \Omega^{1}(M) \to \Omega^{0}(M) \oplus \Omega^{1}(M), \;\:
(f, \beta) \mapsto (-\inner{\alpha,\beta}, \:
f \alpha + *(\alpha \wedge \beta \wedge *\varphi)) ,
\end{equation*}
and the Dirac operator is identified with
\begin{equation}
\label{diracmodeleq}
\Omega^{0}(M) \oplus \Omega^{1}(M) \to \Omega^{0}(M) \oplus \Omega^{1}(M), \;\:
(f, \beta) \mapsto (d^{*}\beta, \, df + *(d\beta \wedge *\varphi)) .
\end{equation}
If $M$ is a \gmfd{} we can see directly from (\ref{diracmodeleq}) that
$\dirac^{2}$ is identified
with the Hodge Laplacian on $\Omega^{0}(M) \oplus \Omega^{1}(M)$ (this could
also be proved in the vein of proposition \ref{hodgelapprop}, using
a Weitzenböck formula). This allows us to deduce the index of the Dirac
operator from that of the Laplacian (proposition \ref{lapindprop}).

Let $\harm^{S}_{0}$ be the bounded harmonic spinors on the \gmfd{} $M^{7}$,
and $\harm^{S}_{\infty}$ the translation-invariant harmonic spinors on
the cylinder $X \times \bbr$. We can consider $\rho \harm^{S}_{\infty}$
as a space of spinors on $M$ supported on the cylindrical part, where
$\rho$ is a cut-off function for the cylinder.

\begin{prop}
\label{diracprop}
Let $M$ be an EAC \gmfd{} with rate $\delta_{0}$, $k \geq 1$ and
$0 < \delta < \min\{\epsilon_{1}, \delta_{0}\}$.
\begin{subequations}
\begin{equation}
\label{topdiraceq}
\dirac : \holdad{k+1}(S) \to \holdad{k}(S)
\end{equation}
is injective and its image is the \ltwoorth{} complement of $\harm^{S}_{0}$,
while
\begin{gather}
\label{middiraceq}
\dirac : \holdad{k+1}(S) \oplus \rho \harm^{S}_{\infty} \to \holdad{k}(S)
\end{gather}
\end{subequations}
is surjective with kernel $\harm^{S}_{0}$.
\end{prop}

\begin{proof}
The argument is analogous to that on page \pageref{padeq}.
$\harm^{S}_{\infty} \cong \harm^{0}_{\infty} \oplus \harm^{1}_{\infty} \cong
2\harm^{0}_{X} \oplus \harm^{1}_{X}$, so has even dimension (as $X$ is Kähler).
Let $2i = \dim \harm^{S}_{\infty}$.
By proposition \ref{lapindprop} 
\[ \triangle : \holdad{k+1}(S) \to \holdad{k}(S) \]
has index $-2i$, so the index of (\ref{topdiraceq}) is $-i$.
Hence (\ref{middiraceq}) 
has index $+i$.
The kernel of (\ref{middiraceq}) is contained in $\harm^{S}_{0}$,
while by integration by parts the image of (\ref{topdiraceq}) is contained
in the \ltwoorth{} complement of $\harm^{S}_{0}$. Hence
\begin{equation}
\label{indexsqueezeeq}
{\setlength\arraycolsep{2pt}
\begin{array}{rcl}
-i &\leq &-\dim \harm^{S}_{0} , \\
i &\leq &\dim \harm^{S}_{0} .
\end{array}}
\end{equation}
Since equality holds the result follows.
\end{proof}

\begin{rmk}
If $M$ has a single end then $\dim \harm^{S}_{\infty} = 2 + b^{1}(X)$,
while corollary \ref{degonehodgecor} implies that
$\dim \harm^{S}_{0} = 1 + b^{1}(M)$.
Hence the equality in (\ref{indexsqueezeeq})
can be seen as a consequence of proposition \ref{realsubprop}.
\end{rmk}

\subsection{The slice}

Fix $k \geq 1$, $\delta > 0$, $\alpha \in (0,1)$ and $\varphi \in \calxa$.
We find a direct complement $\slt$ in $T_{\varphi}\ztreq$ to the
tangent space of the $\calq{D}$-orbit.
Then we define a submanifold $\mathcal{S} \subseteq \ztreq$
whose tangent space at $\varphi$ is $\slt$. We will use $\mathcal{S}$ as a
slice in $\ztreq$ for the $\calq{D}$-action at $\varphi$.

The fixed torsion-free \gstr{} $\varphi$ is used to define an EAC metric and a
Hodge star.
It also defines type decompositions of the exterior bundles and spaces of
harmonic forms, as described in subsection \ref{harmholsub}.
The relevant decompositions of the exterior powers of the cotangent bundle are
\begin{align*}
\Lambda^{2} &= \Lambda^{2}_{7}
\oplus \Lambda^{2}_{14} \\
\Lambda^{3} &= \Lambda^{3}_{1} \oplus
\Lambda^{3}_{7} \oplus \Lambda^{3}_{27} \\
\Lambda^{4} &= \Lambda^{4}_{1} \oplus
\Lambda^{4}_{7} \oplus \Lambda^{4}_{27} \\
\Lambda^{5} &= \Lambda^{5}_{7} \oplus
\Lambda^{5}_{14} ,
\end{align*}
where $\Lambda^{m}_{d}$ is a subbundle of the exterior cotangent bundle
$\Lambda^{m}$ of rank $d$. Its sections $\Omega^{m}_{d}(M)$ are the
`type $d$ $m$-forms'. Note that the map
$TM \to \Lambda^{2}_{7}, \; v \mapsto v \lrcorner \varphi$ is a bundle
isomorphism. Therefore Lie derivatives
$\mathcal{L}_{V}\varphi = d(V \lrcorner \varphi)$ are precisely exterior
derivatives of $2$-forms of type $7$.

Having restricted our attention to \gstr s in $\ztreq$ is convenient
since the asymptotic values of elements of $T_{\varphi}\ztreq$ are harmonic.
Recall the notation for harmonic forms from subsection \ref{hodgesub},
in particular that $\harm^{m}_{0}$
and $\harm^{m}_{\infty}$ denote the spaces of bounded harmonic forms on
$M$ and harmonic translation-invariant forms on $X \times \bbr$ respectively.

For convenience we identify translation-invariant $3$-forms on $X \times \bbr$
with pairs of $3$- and $2$-forms on $X$ by
$\sigma + dt \wedge \tau \leftrightarrow (\sigma, \tau)$. This identifies the
tangent spaces $\harm_{SU}$ and $\harm_{SU,A}$ of $\mathcal{Q}$ and
$\mathcal{Q}_{A}$ with subspaces $\harm^{3}_{SU}, \harm^{3}_{SU,A} \subseteq
\harm^{3}_{\infty}$.
Let
\begin{equation*}
\ztre_{cyl} \subseteq \holdad{k}(\Lambda^{3}) \oplus \rho \harm^{3}_{SU,A}
\end{equation*}
be the subspace of closed forms. Clearly
$T_{\varphi}\ztreq \subseteq \ztre_{cyl}$, and we show below that
equality holds.
The tangent space to the pre-moduli space of torsion-free \gstr s at $\varphi$
will turn out to be the subspace
\begin{equation*}
\label{hcyleq}
\harm^{3}_{cyl} \subseteq \ztre_{cyl}
\end{equation*}
of harmonic forms. This is exactly the subspace of elements of 
$\harm^{3}_{0}$ which are tangent to cylindrical deformations of the \gstr{},
i.e. whose boundary values lie in $\harm^{3}_{SU}$. The boundary map
$B : \harm^{3}_{0} \to \harm^{3}_{\infty}$ maps $\harm^{3}_{cyl}$ precisely
onto $\harm^{3}_{SU,A}$. Together with the Hodge decomposition
(\ref{hodgedecompeq}) it follows that
\begin{equation}
\label{zcylspliteq}
\ztre_{cyl} = \harm^{3}_{cyl} \oplus \holdad{k}[d\Lambda^{2}] .
\end{equation}

\begin{rmk}
$d\holdad{k+1}(\Lambda^{m})$ is the space of exterior derivatives of decaying
forms, while we use $\holdad{k}[d\Lambda^{m}]$ to denote the space of exact
decaying forms.
$d\holdad{k+1}(\Lambda^{m}) \subseteq \holdad{k}[d\Lambda^{m}]$ is a closed
subspace of finite codimension.
\end{rmk}

\begin{lem}
\label{zqmfdlem}
$\ztre_{\mathcal{Q}}$ is a manifold, and $T_{\varphi}\ztreq = \ztre_{cyl}$.
\end{lem}

\begin{proof}
If $\psi$ is a $3$-form asymptotic to an element
$(\Omega',\omega') \in \mathcal{Q}_{A}$
then the condition $[\Omega'] \in A^{3}$ implies that
$d\psi \in d\holdad{k}(\Lambda^{3})$. Therefore
\[ d : \holdad{k}(\Lambda^{3}) + \rho \mathcal{Q}_{A} \to
d\holdad{k}(\Lambda^{3}) \]
is a submersion, and the result follows from the implicit function theorem.
\end{proof}

Let $\dkaq$ be the group of diffeomorphisms of $M$ which are isotopic to the
identity, and $\holdad{k+1}$-asymptotic to a cylindrical automorphism of the
cylindrical \gstr{} $\Omega + dt \wedge \omega$. The elements of a
neighbourhood of the identity in $\dkaq$ can be written as
$\exp (V + \rho V_{\infty})$, where $V$ is a $\holdad{k+1}$ vector field on $M$
and $V_{\infty}$ is a translation-invariant vector field on $X \times \bbr$
with $\mathcal{L}_{V_{\infty}}(\Omega + dt \wedge \omega) = 0$, i.e.
$V_{\infty} \in (\harm^{1}_{\infty})^{\sharp}$.
Therefore if we let
\begin{equation*}
\label{deq}
D = \rho (\harm^{1}_{\infty})^{\sharp} \lrcorner \varphi
\subseteq \Omega^{2}_{7}(M)
\end{equation*}
then the tangent space to the $\dkaq$-orbit at $\varphi$ is
\[ \{ \mathcal{L}_{V + \rho V_{\infty}}\varphi : 
V \in \holdad{k+1}(TM), V_{\infty} \in (\harm^{1}_{\infty})^{\sharp} \} =
d(\holdad{k+1}(\Lambda^{2}_{7}) \oplus D) . \]
We take the tangent space $\slt$ to the slice at $\varphi$ to be the kernel
of the formal adjoint of $d : \Omega^{2}_{7}(M) \to \Omega^{3}(M)$, i.e.

\begin{defn}
Let $\slt$ be $\ker \pi_{7}d^{*}$ in $\ztre_{cyl}$.
\end{defn}

\begin{lem}
\label{bryantlem}
If $\beta \in \Omega^{3}_{27}(M)$ then $\pi_{7}d\beta = 0$ if and only if
$\pi_{7}d^{*}\beta = 0$.
\end{lem}

\begin{proof}
An instance of \cite[Proposition $3$]{bryant03}.
\end{proof}

\begin{defn}
Let $W = \holdad{k}[d\Lambda^{2}] \cap \Omega^{27}(M)$.
\end{defn}

\begin{prop}
\label{sltprop}
$\slt = \harm^{3}_{cyl} \oplus W$, and
\begin{equation*}
\ztre_{cyl} = d(\holdad{k+1}(\Lambda^{2}_{7}) \oplus D) \oplus \slt .
\end{equation*}
\end{prop}

\begin{proof}
$\slt$ obviously contains $\harm^{3}_{cyl}$, and it also contains $W$ by lemma
\ref{bryantlem}. By integration by parts $\slt$ is \ltwoorth{} to
$d(\holdad{k+1}(\Lambda^{2}_{7}) \oplus D)$. Hence by (\ref{zcylspliteq})
it suffices to show
\begin{equation}
\label{wspliteq}
\holdad{k}[d\Lambda^{2}] = d(\holdad{k+1}(\Lambda^{2}_{7}) \oplus D) \oplus W .
\end{equation}
We can identify the spinor bundle $S$ with $\Lambda^{0} \oplus \Lambda^{2}_{7}$
and with $\Lambda^{3}_{1\oplus7}$
(shorthand for $\Lambda^{3}_{1} \oplus \Lambda^{3}_{7}$)
so that the Dirac operator
$\dirac : \Gamma(S) \to \Gamma(S)$ is identified with
\begin{equation*}
\Omega^{0}(M) \oplus \Omega^{2}_{7}(M) \to \Omega^{3}_{1 \oplus 7}(M), \;\:
(f, \eta) \mapsto \pi_{1\oplus7}d\eta + *(df \wedge \varphi) .
\end{equation*}
If $\beta \in \holdad{k}[d\Lambda^{2}]$ then by surjectivity of the Dirac
operator map (\ref{middiraceq})
\[ \pi_{1\oplus7}\beta = \pi_{1\oplus7}d\eta + *(df \wedge \varphi) \]
for some $\eta \in \holdad{k+1}(\Lambda^{2}_{7}) \oplus D$,
$f \in \holdad{k+1}(\Lambda^{0}) \oplus \rho\harm^{0}_{\infty}$.
By integration by parts $df = 0$. Hence
$\beta - d\eta$ is exact and of type $27$, i.e. $\beta - d\eta \in W$.
\end{proof}

We want to take as our slice for the $\cald{D}$-action at $\varphi$ a
submanifold $\mathcal{S} \subseteq \ztreq$ with
$T_{\varphi}\mathcal{S} = \slt$.
To this end we pick a map
\begin{equation}
\label{expdefeq}
\exp : U \to \ztreq
\end{equation}
on a neighbourhood $U$ of $0$ in $\ztre_{cyl} = T_{\varphi}\ztreq$,
such that $D\exp_{0} = id$. We also insist that $\exp$ maps
decaying forms to decaying forms, and smooth forms to smooth forms.
We can do this since by (\ref{zcylspliteq}) the decaying forms have a
finite-dimensional complement of smooth forms in $\ztre_{cyl}$.
We then choose
\begin{equation}
\label{slicedefeq}
\mathcal{S} = \exp (\slt \cap U) .
\end{equation}

\subsection{The pre-moduli space is a manifold}

Let $\cald{R} \subseteq \mathcal{S}$ be the subset of $\holdad{k}$
torsion-free \gstr s.
$\cald{R}$ is the pre-moduli space of torsion-free \gstr s near $\varphi$.
In order to show that $\cald{R}$ is a submanifold we will
exhibit it as the zero set of a function $F$ with surjective derivative, and
apply the implicit function theorem.

Recall that by theorem \ref{graythm} a \gstr{} $\psi$ is torsion-free if and
only if $d\psi = 0$ and $d^{*}_{\psi}\psi = 0$. To emphasise the non-linearity
of the second condition we define

\begin{defn}
For a \gstr{} $\psi$ on $M$ let $\Theta(\psi) = *_{\psi}\psi$.
\end{defn}

So with this notation
\[ \cald{R} = \{ \psi \in \mathcal{S} : d\Theta(\psi) = 0 \} . \]
If $\psi \in \ztreq$ then $\psi$ is asymptotic to a torsion-free
cylindrical \gstr, so $d\Theta(\psi)$ decays. Moreover, by definition elements
of $\ztreq$ are asymptotic to elements of
$\mathcal{Q}_{A} \subseteq \mathcal{Q}$.
Therefore $\Theta(\psi)$ is asymptotic
to $\half \omega'^{2} - dt \wedge \hat\Omega'$, with $[\omega'^{2}] \in A^{4}$
(cf. (\ref{qaeq}) and remark \ref{cyldualrmk}).
This implies that $d\Theta(\psi) \in d\holdad{k+1}(\Lambda^{4})$.

$\Theta : \Lambda^{3}_{+}T^{*}M \to \Lambda^{4}T^{*}M$ is point-wise smooth,
so by the chain rule
\begin{equation*}
\ztreq \to d\holdad{k+1}(\Lambda^{4}), \;\: \psi \to d\Theta(\psi)
\end{equation*}
is a smooth function.
We need to adjust this map to obtain a function with surjective derivative.
If $\beta$ is a $3$-form such that
$d^{*}\beta \in d^{*}\holdad{k+1}(\Lambda^{3})$ then by Hodge decomposition
(\ref{hodgedecompeq}) there is a unique
$\beta_{E} \in \holdad{k+1}[d\Lambda^{2}]$ such that
$d^{*}\beta = d^{*}\beta_{E}$.
We can think of $\beta_{E}$ as the exact part of $\beta$.
The image of $\beta_{E}$ under the projection
$\holdad{k+1}[d\Lambda^{2}] \to W$ in the direct sum decomposition
(\ref{wspliteq}) is the unique $P(\beta) \in W$ such that
\[ d^{*}\beta - d^{*}P(\beta) \in
d^{*}d(\holdad{k+1}(\Lambda^{2}_{7}) \oplus D) . \]

\begin{defn}
For $\psi$ close to $\varphi$ in $\ztreq$ let
\begin{equation}
F(\psi) = P(*\Theta(\psi)) .
\end{equation}
\end{defn}

$F(\psi)$ is essentially a component of the exact part of $*\Theta(\psi)$,
so $d\Theta(\psi) = 0 \Rightarrow F(\psi) = 0$. In order to show that
the converse also holds, so that we do not `lose any information' by
considering zeros of $F$ instead of $\psi \mapsto d\Theta(\psi)$, we need
a simple algebraic lemma.

\begin{lem}
\label{grayiilem}
Let $\psi$ a \gstr{} on $M$ with $d\psi = 0$. Then for any vector
field $V$
\begin{equation*}
d\Theta(\psi) \wedge (V \lrcorner \psi) = 0 .
\end{equation*}
\end{lem}

\begin{proof}
By \cite[Lemma $11.5$]{salamon89} there is for any \gstr{} $\psi$ a linear
relation between $\pi_{7}d\psi$ and $\pi_{7}d^{*}\psi$
(where the type component and codifferential are defined by $\psi$).
This implies the result.
\end{proof}

\begin{prop}
\label{zerosetprop}
For $\psi \in \ztreq$ sufficiently close to $\varphi$, $\psi$ is
torsion-free if and only if $F(\psi) = 0$.
\end{prop}

\begin{proof}
By definition $F(\psi) = 0$ implies that
\[ *d\Theta(\psi) = d^{*}d\eta \]
for some $\eta \in \holdad{k+1}(\Lambda^{2}_{7}) \oplus D$.
We need to show that $d^{*}d\eta = 0$.
By integration by parts it would suffice to show that $\pi_{7}d^{*}d\eta = 0$,
but unfortunately there is no \textit{a priori} reason why that should
be the case. Nevertheless, if $\psi$ is close to $\varphi$ then
$\norm{\pi_{7}d^{*}d\eta}$ is small relative to $\norm{d^{*}d\eta}$,
for if $V \lrcorner \varphi \in \Lambda^{2}_{7}$ then 
\begin{multline*}
\inner{d^{*}d\eta, V \lrcorner \varphi} \; = \;
\inner{d^{*}d\eta, V \lrcorner (\varphi - \psi)} \; \leq 
\norm{d^{*}d\eta}\norm{\psi - \varphi}\norm{V} = 
\third \norm{d^{*}d\eta}\norm{\psi - \varphi}\norm{V \lrcorner \varphi}
\end{multline*}
by lemma~\ref{grayiilem}.
Hence point-wise
\begin{equation}
\label{pwboundeq}
\norm{\pi_{7}d^{*}d\eta} \leq \third \norm{\psi - \varphi}\norm{d^{*}d\eta} .
\end{equation}

On the other hand we can find a reverse inequality for weighted
$L^{2}$-norms: if we pick $0 < \delta' < \delta$ then
there is a constant $C > 0$ such that
\begin{equation}
\label{sobboundeq}
\norm{\pi_{7}d^{*}d\eta}_{\sob{\delta'}} \geq
C \norm{d^{*}d\eta}_{\sob{\delta'}}
\end{equation}
for any $\eta \in \sobx{2}{\delta'}(\Lambda^{2}_{7}) \oplus D$.
To prove this inequality we use weighted Sobolev $\sobx{2}{\delta'}$ norms,
defined analogously to the weighted Hölder norms in definition \ref{wholddef}.
$\holdad{k+1} \hookrightarrow \sobx{2}{\delta'}$ continuously for $k \geq 1$.
Since $\pi_{7}d^{*}d : \Omega^{2}_{7}(M) \to \Omega^{2}_{7}(M)$
is elliptic \cite[Theorem $6.2$]{lockhart85} (discussed on page
\pageref{tiopeq}) ensures that
\[ \pi_{7}d^{*}d : \sobx{2}{\delta'}(\Lambda^{2}_{7}) \oplus D
\to \sob{\delta'}(\Lambda^{2}_{7}) \]
is Fredholm. Therefore for a suitable choice of $C$ the inequality
(\ref{sobboundeq}) holds for $\eta$ in a
complement of the kernel of $\pi_{7}d^{*}d$. By integration by parts the
elements of the kernel are closed, so (\ref{sobboundeq}) holds
for all $\eta \in \sobx{2}{\delta'}(\Lambda^{2}_{7}) \oplus D$.
Combining the inequalities (\ref{pwboundeq}) and (\ref{sobboundeq}) gives that
if $\norm{\psi - \varphi}_{C^{0}} < 3C$ then $d^{*}d\eta = 0$. Hence
\[ F(\psi) = 0 \Rightarrow d\Theta(\psi) = 0 . \qedhere \]
\end{proof}

Next we show that $F : \mathcal{S} \to W$ satisfies the hypotheses of the
implicit function theorem.

\begin{prop}
$F : \ztreq \to W$ is smooth with derivative
\begin{equation}
\label{derfeq}
DF_{\varphi} : \ztre_{cyl} \to W,
\;\: \beta \mapsto P(\fourthird\pi_{1}\beta + \pi_{7}\beta - \pi_{27}\beta) .
\end{equation}
\end{prop}

\begin{proof}
According to \cite[Lemma $3.1.1$]{joyce96-I} the derivative at $\varphi_{0}$
of the point-wise model
$\Theta : \Lambda^{3}_{+}(\bbr^{7})^{*} \to \Lambda^{4}(\bbr^{7})^{*}$
is
\begin{equation}
\label{derthetaeq}
D\Theta_{\varphi_{0}} :
\Lambda^{3}(\bbr^{7})^{*} \to \Lambda^{4}(\bbr^{7})^{*}, \;\:
\beta \mapsto *\fourthird\pi_{1}\beta + *\pi_{7}\beta - *\pi_{27}\beta ,
\end{equation}
and the result follows by the chain rule.
\end{proof}

\begin{prop}
\label{cylfderprop}
$DF_{\varphi} : \ztre_{cyl} \to W$ is $0$ on
$d(\holdad{k+1}(\Lambda^{2}_{7}) \oplus D) \oplus \harm^{3}_{cyl}$
and $-id$ on~$W$.
\end{prop}

\begin{proof}
For any $V \lrcorner \varphi \in
\holdad{k+1}(\Lambda^{2}_{7}) \oplus D$
\[ d^{*}\!*\!D\Theta_{\varphi}(d(V \lrcorner \varphi)) =
*d(D\Theta_{\varphi})(\mathcal{L}_{V}\varphi) = 
*\mathcal{L}_{V}(d\Theta(\varphi)) = 0\]
so $DF_{\varphi}(d(V \lrcorner \varphi)) = 0$.
The type components of harmonic forms are harmonic, so
for any $\chi \in \harm^{3}_{cyl}$ we have
$d^{*}(\frac{4}{3}\pi_{1}\chi + \pi_{7}\chi - \pi_{27}\chi) = 0$,
and hence $DF_{\varphi}(\chi) = 0$.
If $\kappa \in W$ then $DF_{\varphi}(\kappa) = P(-\kappa) = -\kappa$,
as $P$ is a projection to $W$.
\end{proof}

We have taken the pre-moduli space $\cald{R}$ near $\varphi$ to consist of the
torsion-free EAC \gstr s in the slice $\mathcal{S}$. We can, however, only
prove that it has the properties we want close to $\varphi$. We will therefore
repeatedly replace $\mathcal{S}$ by a neighbourhood of $\varphi$ in
$\mathcal{S}$ in order to ensure that $\cald{R} \subseteq{S}$ has the
desired properties.
The first step is to ensure that $\cald{R}$ is a manifold.
Proposition \ref{zerosetprop} shows that if $\mathcal{S}$ is sufficiently small
then $\cald{R}$ is the zero set of $F$ in $\mathcal{S}$.
Applying the implicit function theorem to $F : \mathcal{S} \to W$ we deduce

\begin{prop}
\label{premodprop}
If the slice $\mathcal{S}$ near $\varphi$ is shrunk sufficiently then
the space $\cald{R}$ of torsion-free EAC \gstr s in $\mathcal{S}$ is a manifold.
Its tangent space at $\varphi$ is $\harm^{3}_{cyl}$.
\end{prop}

\subsection{Regularity}
\label{cylregsub}

In order to use the pre-moduli space $\cald{R}$ to construct a moduli
space of \emph{smooth} \gstr s we require two regularity results.
We show that if the slice $\mathcal{S}$ is sufficiently small
then the elements of the pre-moduli space $\cald{R} \subseteq \mathcal{S}$,
which \textit{a priori} are merely $\holdad{k}$, are smooth and EAC.
We also need to show that isometries of EAC metrics are EAC.
To prove the regularity of elements of $\cald{R}$
we use a regularity result about solutions of elliptic operators which are
asymptotically translation-invariant in a weighted Hölder sense.

\begin{defn}
A differential operator $A$ on $M$ is \emph{\holdadasy{l}}
to the translation-invariant operator $A_{\infty}$ on the cylinder
$X \times \bbr$ if the difference between the coefficients
of $A$ and $A_{\infty}$ (in the sense of (\ref{asycoeffeq})) is~$\holdad{l}$.
\end{defn}

The theorem below is similar to e.g. \cite[Theorem $6.3$]{mazya78},
and can be obtained from interior estimates like those of Morrey 
\cite[Theorems $6.2.5$ and $6.2.6$]{morrey66}.

\begin{thm}
\label{cylregthm}
Let $M^{n}$ be an asymptotically cylindrical manifold, $E$ a vector bundle
on $M$ associated to $TM$,
and $A$ a linear elliptic order $r$ differential operator on the sections
of $E$ that is $\holdad{l}$-asymptotic to a translation-invariant operator
with $\holda{l}$ coefficients.
If $s \in \holdad{r}(E)$ and $As \in \holdad{l}(E)$ then
$s \in \holdad{l+r}(E)$, and there is a constant $C > 0$ independent of $s$
such that
\begin{equation*}
\norm{s}_{\holdad{l+r}} < C\left( \norm{As}_{\holdad{l}} +
\norm{s}_{C^{0}_{\delta}} \right) .
\end{equation*}
\end{thm}

Now consider again a \gmfd{} $M$ with a torsion-free \gstr{}~$\varphi$,
and the pre-moduli space $\cald{R}$ of torsion-free \gstr s
in the slice $\mathcal{S} = \exp(\slt \cap U) \subseteq \ztreq$ for the
$\calq{D}$-action at $\varphi$.
We use theorem \ref{cylregthm} in a boot-strapping argument to show
that the elements of $\cald{R}$ are EAC. A priori they are
$\holdad{k}$-asymptotic to elements of $\mathcal{Q}_{A}$.
Like in proposition \ref{premodprop} we can only work close to $\varphi$, and
must replace $\mathcal{S}$ by a neighbourhood of $\varphi$ in $\mathcal{S}$.

\hfuzz=5pt

\begin{prop}
\label{cylregprop}
If the slice $\mathcal{S}$ near $\varphi$ is shrunk sufficiently then the
pre-moduli space $\cald{R} \subseteq \mathcal{S}$ consists of
smooth exponentially asymptotically translation-invariant forms.
\end{prop}

\hfuzz=0pt

\begin{proof}
We want to show that if $\psi \in \mathcal{S}$ is sufficiently close to
$\varphi$ and $d\Theta(\psi) = 0$ then $\psi$ is smooth and
exponentially asymptotically translation-invariant.
Write $\psi = \varphi + \beta$.
\[ D(d\!*\!d\Theta)_{\varphi} =
-dd^{*} \circ (\fourthird\pi_{1} + \pi_{7} - \pi_{27}) , \]
so we can write
\[ d\!*\!d\Theta(\varphi + \beta) =
-dd^{*}(\fourthird\pi_{1}\beta + \pi_{7}\beta - \pi_{27} \beta) +
P(\beta, \nabla \beta, \nabla^{2}\beta) + R(\beta, \nabla \beta) , \]
where $P$ consists of the quadratic terms of $d\!*\!d\Theta(\varphi+\beta)$
that involve $\nabla^{2}\beta$, and $R$ consists of the remaining quadratic
terms. $P$ and $R$ depend only point-wise on their arguments, and $P$
is linear in $\nabla^{2}\beta$.

By the definition of the map $\exp$ (\ref{expdefeq}) we can write
$\beta = \kappa + \gamma$, with $\kappa \in W$, and
$\gamma$ smooth and exponentially asymptotic to some element of
$\mathcal{Q}_{A}$.
As $\kappa$ is closed of type $27$
\[ -dd^{*}(\fourthird\pi_{1}\kappa + \pi_{7}\kappa - \pi_{27}\kappa) =
\triangle \kappa . \]

Considering $\beta$ and $\nabla \beta$ as fixed we can define a
second-order linear differential operator
$A : \zeta \mapsto P(\beta, \nabla\beta, \nabla^{2}\zeta)$. Then
the condition $d\!*\!d\Theta(\psi) = 0$ is equivalent to 
\begin{equation}
\label{lineareq2}
(\triangle + A)\kappa = -R +
dd^{*}(\fourthird\pi_{1}\gamma + \pi_{7}\gamma - \pi_{27}\gamma) - A\gamma .
\end{equation}

If $\beta = 0$ then $A = 0$, so $\triangle + A$ is elliptic. Ellipticity
is an open condition, so $\triangle + A$ is in fact elliptic for any
$\beta$ sufficiently small in the uniform norm.

Now suppose $\kappa$ is $\holdad{l+1}$ and is a solution of
(\ref{lineareq2}). $R$ and the coefficients of $A$ depend smoothly on $\kappa$
and $\nabla \kappa$. Therefore $\triangle + A$ and the RHS of (\ref{lineareq2})
are $\holdad{l}$-asymptotically translation-invariant.
Since the RHS of (\ref{lineareq2}) is decaying \textit{a priori} it is
$\holdad{l}$.
If $\beta$ is sufficiently small that $\triangle + A$ is elliptic then by
theorem \ref{cylregthm} $\kappa$ is $\holdad{l+2}$.
Since $\kappa$ is $\holdad{1}$ it is $\holdad{l}$ for all
$l$ by induction.

Hence $\psi = \varphi + \kappa + \gamma$ is smooth and exponentially
asymptotically translation-invariant.
\end{proof}

Myers and Steenrod \cite{myers39} show that 
any isometry of smooth Riemannian manifolds is smooth.
We wish to generalise this, and show that isometries of EAC metrics are EAC
(in the sense of definition \ref{cyldiffdef}).

We can think of the choice of diffeomorphism $M_{\infty} \to X \times \bbrp$
in definition \ref{cylenddef} as defining a `cylindrical-end structure',
and of two such structures as being `$\delta$-equivalent'
if they differ by a rate~$\delta$ EAC diffeomorphism -- then they define
equivalent notions of exponential translation-invariance etc.
The next proposition can be interpreted as stating that
the cylindrical-end structure of an EAC manifold can be recovered from the
metric.

\begin{prop}
\label{eacisomprop}
Any isometry of EAC (rate $\delta > 0$) manifolds is $C^{\infty}$
and EAC with rate~$\delta$.
\end{prop}

\begin{proof}[Sketch proof.]
By \cite[Theorem 8]{myers39} the isometries are $C^{\infty}$,
so we just need to prove that they are also EAC.

Let $M$ be a manifold with a Riemannian metric $g$. We need to show that
if for $i = 1$, $2$ $M_{i,\infty} \subseteq M$ have compact complements,
and $\Psi_{i} : M_{i,\infty} \to X_{i} \times \bbrp$ are diffeomorphisms
defining cylindrical-end structures with respect to which $g$ is EAC, then
$\Psi_{1} \circ \Psi_{2}^{-1}$ is EAC.

Consider the space $R$ of half-lines in $M$, i.e. equivalence classes of unit
speed globally distance-minimising geodesic rays $\gamma : [0,\infty) \to M$,
where two rays are equivalent if one is a subset of the other.
We can define a distance function on $R$ by
\begin{equation}
d([\gamma],[\sigma]) = \lim_{u \to \infty} \inf_{v} d(\gamma(u), \sigma(v)) .
\end{equation}
$g$ is pushed forward to an EAC metric on $X_{i} \times \bbrp$ by $\Psi_{i}$.
It is straight-forward to solve the geodesic equation in local coordinates
to show that for each $x \in X_{i}$ there is a unique half-line $[\gamma]$
such that the $X$-component of $\gamma(u)$ approaches $x$ as $u \to \infty$.
$\Psi_{i}$ therefore induce isometries $\Xi_{i} : R \to X_{i}$.
Then $\Xi_{1} \circ \Xi_{2}^{-1}$ is smooth by \cite[Theorem 8]{myers39}.
If $t_{i}$ is the $\bbrp$-coordinate on $X_{i} \times \bbrp$ then
$grad (t_{i}) - \frac{d\gamma}{du}$ decays exponentially as $u \to \infty$ for
each half-line $[\gamma]$, so $\Psi_{1} \circ \Psi_{2}^{-1}$
is exponentially asymptotic to
$(x,t) \mapsto (\Xi_{1} \circ \Xi_{2}^{-1}(x), t+h)$ for some $h \in \bbr$.
\end{proof}

\subsection{Constructing the moduli space}
\label{modconstrsub}

For each $\varphi \in \calxa$ we have constructed a pre-moduli
space~$\cald{R}$. $\cald{R}$ is a manifold, its elements are smooth and EAC,
and its tangent space at $\varphi$ is~$\harm^{3}_{cyl}$.
We now want to use slice arguments to show that we can take the pre-moduli
spaces $\cald{R}$ as coordinate charts to define a smooth structure on
$\cald{M}$.
In \cite{ebin70} Ebin gives a very detailed account of a slice construction
on a compact manifold. While the basic idea of the slice is the same,
it is not so convenient for our purposes to attempt to extend Ebin's
arguments to the EAC setting. It is much easier to study the charts for
$\calp{M}$ in terms of the projection to the de Rham cohomology
which appears in the statement of the main theorem \ref{maincylthm}.
\begin{equation*}
\pi_{\mathcal{M}} : \calp{X} \to H^{3}(M) \times H^{2}(X), \;\:
\varphi \mapsto ([\varphi], [B_{e}(\varphi)]) ,
\end{equation*}
where we use $B_{a}(\beta) + dt \wedge B_{e}(\beta)$ to denote the
asymptotic limit of an asymptotically translation-invariant form $\beta$.

We first check that $\pi_{\mathcal{M}}$ is an embedding on $\cald{R}$.
If we allow ourselves to shrink $\cald{R}$ this
amounts to showing that the derivative of $\pi_{\mathcal{M}}$ at $\varphi$ is
injective, and from the definition the derivative is
plainly
\begin{equation*}
(\pi_{H}, \pi_{H,e}) : \harm^{3}_{cyl} \to H^{3}(M) \times H^{2}(X), \;\:
\beta \mapsto ([\beta], [B_{e}(\beta)]) .
\end{equation*}
The kernel consists of harmonic, exact, decaying forms, so
it is trivial.

Recall from subsection \ref{outlinesub} that we chose a pre-moduli space
$\mathcal{Q}$ near the Calabi-Yau structure $(\Omega,\omega)$ on $X$ defined
by the asymptotic limit of $\varphi$, and that $\calq{X} \subseteq \calxa$
is the subset of smooth torsion-free \gstr s whose asymptotic limits lie in
$\mathcal{Q}$. $\calq{D} \subseteq \cald{D}$ is the subgroup of smooth EAC
diffeomorphisms of $M$ whose asymptotic limits lie in the automorphism group
of $(\Omega,\omega)$. $\calq{D}$~acts on $\calq{X}$ by proposition
\ref{autinvprop}, and as an intermediate step for our slice result we prove
that $\cald{R}$ is a coordinate chart for $\calq{X}/\calq{D}$.

\begin{prop}
\label{qsliceprop}
The natural map $\cald{R} \to \calq{X}/\calq{D}$ is a homeomorphism onto
a neighbourhood of $\varphi\calq{D}$.
\end{prop}

\begin{proof}
Let $\mathcal{X}_{k,\mathcal{Q}}$ be the space of torsion-free \gstr s in
$\ztreq$.
By proposition \ref{cylfderprop} a small neighbourhood of $\varphi$ in
$\mathcal{X}_{k,\mathcal{Q}}$ is a manifold, and its tangent space at $\varphi$
is $d(\holdad{k+1}(\Lambda^{2}_{7}) \oplus D) \oplus \harm^{3}_{cyl}$.
The first term is the tangent space to the
$\dkaq$-orbit of $\varphi$, so the derivative at $(\varphi, id)$ of
\begin{equation*}
\cald{R} \times \dkaq \to \ztreq, \;\: (\beta, \phi) \mapsto \phi^{*}\beta
\end{equation*}
is surjective. By the submersion theorem it is an open map on a neighbourhood
of $(\varphi,id)$.

Thus if $U \subseteq \cald{R}$ is a small neighbourhood of $\varphi$ then
$U\dkaq$ is open in $\mathcal{X}_{k,\mathcal{Q}}$.
If $\psi \in \cald{R}$, $\phi \in \dkaq$ and $\phi^{*}\psi$ is smooth
and EAC then $\phi$ is also smooth and EAC by proposition \ref{eacisomprop}.
Therefore $U\dkaq \cap \calq{X} = U\calq{D}$, and $U\calq{D}$ is open
in~$\calq{X}$.

Hence $\cald{R} \to \calq{X}/\calq{D}$ is an open map. It is also injective,
since $\pi_{\mathcal{M}}$ is $\cald{D}$-invariant and
injective on $\cald{R}$.
\end{proof}

For our argument to work we may need to shrink $\mathcal{Q}$ by replacing it
with a neighbourhood of $(\Omega,\omega)$ in $\mathcal{Q}$.

\begin{lem}
\label{qmodlem}
If the pre-moduli space $\mathcal{Q}$ of Calabi-Yau structures is shrunk
sufficiently then
$\calq{X}/\calq{D}$ is homeomorphic to a neighbourhood of $\varphi\cald{D}$
in $\calxa/\cald{D}$.
\end{lem}

\begin{proof}
The natural map $f : \calq{X}/\calq{D} \to \calxa/\cald{D}$ is injective
by proposition \ref{autinvprop}.

Let $\mathcal{Y}$ be the space of \cystr s on $X$. The construction of the
moduli space of \cystr s on $X$ uses slice results similar to proposition
\ref{qsliceprop}. This allows us to define on a neighbourhood $U$
of $(\Omega,\omega)$ in $\mathcal{Y}$ a continuous map
\begin{equation*}
P : U \to C^{\infty}(TX) \times \mathcal{Q}, \;\:
(\beta, \gamma) \mapsto (V, \Omega', \omega')
\end{equation*}
such that $(\beta,\gamma) = (\exp V)^{*}(\Omega',\omega')$ for any
$(\beta,\gamma) \in U$.

Let $\mathcal{X}_{U} = \{ \psi \in \calxa : B(\psi) \in U \}$.
If $\psi \in \mathcal{X}_{U}$ let $V = P(B(\psi))$,
$\phi = \exp \rho V \in \cald{D}$ and $g(\psi) = \phi^{*}\psi$.
Then $B(g(\psi)) \in \mathcal{Q}$, so $\psi \in \calq{X}$.
Obviously $f(g(\psi)\calq{D})) = \psi\cald{D}$. Since $f$ is injective
$g$ induces a well-defined map $\mathcal{X}_{U}\cald{D} \to \calq{X}/\calq{D}$.
$g$ is an inverse for $f$ on a neighbourhood of $\varphi\cald{D}$ in
$\calxa/\cald{D}$, so the result follows.
\end{proof}

\begin{thm}
\label{dmodthm}
$\cald{M}$ has a unique smooth structure such that
\[ \pi_{\mathcal{M}} : \cald{M} \to H^{3}(M) \times H^{2}(X) \]
is an immersion.
\end{thm}

\begin{proof}
Proposition \ref{qsliceprop} and lemma \ref{qmodlem} show that
for any $\varphi \in \calxa$ the pre-moduli space $\cald{R}$
is homeomorphic to a neighbourhood $U$ of $\varphi\cald{D}$ in~$\cald{M}$,
and
\[ \pi_{\mathcal{M}} : U \to H^{3}(M) \times H^{2}(X) \]
is a homeomorphism onto an embedded manifold. To prove that we can use the maps
$\cald{R} \to U$ as coordinate charts for $\cald{M}$ we need to check that
the transition functions are smooth. But on an overlap $U_{1} \cap U_{2}$
both charts define the unique smooth structure for which
$\pi_{\mathcal{M}} : U_{1} \cap U_{2} \to H^{3}(M) \times H^{2}(X)$
is an embedding, so we are done.
\end{proof}

If $\delta_{1} > \delta_{2} > 0$ then
$\mathcal{M}_{\delta_{1}} \to \mathcal{M}_{\delta_{2}}$ is injective by
proposition \ref{eacisomprop}, and $\mathcal{M}_{\delta_{1}}$ must be an
open submanifold of $\mathcal{M}_{\delta_{2}}$ since $\pi_{\mathcal{M}}$
is an immersion on both spaces.
Similarly $\cald{M} \to \calp{M}$ is injective for any $\delta > 0$, so
\[ \calp{M} = \bigcup_{\delta > 0} \cald{M} . \]
To finish the proof of the main theorem \ref{maincylthm} it remains only
to observe

\begin{lem}
For any $\delta > 0$ the natural map $\cald{M} \to \calp{M}$ is
a homeomorphism to an open subset.
\end{lem}

\begin{proof}
We need to show that $\cald{M} \to \calp{M}$ is open, i.e. that
if $U \subseteq \calxa$ with $U\cald{D}$ open in $\cald{X}$ then
$U\calp{D}$ is open in $\calp{X}$. By the definition of the topology on
$\calp{X}$ this means that $U\calp{D} \cap \mathcal{X}_{\delta'}$ is open
in $\mathcal{X}_{\delta'}$ for any $\delta' > \delta$.
But proposition \ref{eacisomprop} implies that
$U\calp{D} \cap \mathcal{X}_{\delta'} = U\mathcal{D}_{\delta'}$, which
is open in $\mathcal{X}_{\delta'}$ since $\cald{M} \to \mathcal{M}_{\delta'}$
is a local diffeomorphism.
\end{proof}

This concludes the proof of theorem \ref{maincylthm}.

\subsection{Properties of the moduli space}
\label{propertiessub}

We look at some local properties of the moduli space $\calp{M}$
on an EAC \gmfd{} $M$,
which follow from the existence of a pre-moduli space $\mathcal{R}$
with tangent space $\harm^{3}_{cyl}$.

Firstly, the boundary map $B$ maps $\harm^{3}_{cyl}$ onto $\harm^{3}_{SU,A}$,
so proposition \ref{qasubmfdprop} implies that
$B : \mathcal{R} \to \mathcal{Q}_{A}$ is a submersion.
As $\mathcal{Q}_{A}$ is a homeomorphic to an open set in $\mathcal{N}_{A}$
it follows that $B : \calp{M} \to \mathcal{N}_{A}$ is a submersion onto its
image, and we have proved theorem \ref{modsubmersethm}.

We can now deduce corollary \ref{fixboundarycor}. It suffices to
show that the fibres of $B : \calp{M} \to \mathcal{N}_{A}$ are locally
diffeomorphic to the compactly supported subspace
$H^{3}_{0}(M) \subseteq H^{3}(M)$. 

\begin{lem}
Let $\varphi$ be an EAC torsion-free \gstr{} on $M$, $\mathcal{R}$ the
pre-moduli space of EAC torsion-free \gstr s near $\varphi$ and $\mathcal{Q}$
the pre-moduli space of \cystr s near $B(\varphi)$. The map
\begin{equation*}
\label{fibrecharteq}
\mathcal{R} \to H^{3}(M), \;\: \psi \mapsto [\psi - \varphi]
\end{equation*}
maps a neighbourhood of the fibre of $B :\mathcal{R} \to \mathcal{Q}_{A}$
containing $\varphi$ diffeomorphically to an open subset of $H^{3}_{0}(M)$.
\end{lem}

\begin{proof}
If $\psi$ is in the same fibre as $\varphi$ then $\psi - \varphi$ is
exponentially decaying, so $[\psi - \varphi]$ $\in H^{3}_{0}(M)$.
The tangent space to the fibre at $\varphi$ is the kernel of the derivative
of the submersion $B$, i.e. the subspace $\harm^{3}_{+}$ of decaying forms 
in $\harm^{3}_{cyl} = T_{\varphi}\mathcal{R}$. By theorem
\ref{derhamthm} $\harm^{3}_{+} \cong H^{3}_{0}(M)$, and the result follows.
\end{proof}

Finally, to confirm the formula for the dimension in proposition
\ref{dimmodprop} we just have to compute the dimension of $\harm^{3}_{cyl}$.
Recall from subsection \ref{hodgesub} that $A^{m}$ is the image of
$j^{*} : H^{m}(M) \to H^{m}(X)$,
that $H^{m}(X) = A^{m} \oplus E^{m}$
is an orthogonal direct sum, and
that the harmonic representatives of the summands are denoted by
$\mathcal{A}^{m}$ and $\mathcal{E}^{m}$ respectively.

\hfuzz=4pt
\begin{lem}
\label{hmidlaglem}
Let $M^{4k+3}$ be an oriented EAC 
manifold with cross-section~$X$. Then
$A^{2k+1} \subseteq H^{2k+1}(X)$ has dimension $\half b^{2k+1}(X)$.
\end{lem}
\hfuzz=0pt

\begin{proof}
$H^{2k+1}(X)$ is a symplectic vector space under the Poincar\'e pairing.
In particular $b^{2k+1}(X)$ is even.
$* : H^{2k+1}(X) \to H^{2k+1}(X)$ maps $A^{2k+1}$ isomorphically to its
orthogonal complement $E^{2k+1}$. The Poincar\'e pairing on $H^{2k+1}(X)$
can be expressed as $<\cdot, *\cdot>$, so $A^{2k+1} \subseteq H^{2k+1}(X)$
is a Lagrangian subspace.
\end{proof}

In particular for any EAC \gmfd{} $M$ with cross-section $X$ the long exact
sequence (\ref{relexacteq}) for relative cohomology gives
\begin{equation}
\label{b3cptdimeq}
\dim H^{3}_{0}(M) = b^{3}(M) - \half b^{3}(X) .
\end{equation}

\begin{lem}
\label{dimmodlem}
$\dim \harm^{3}_{cyl} = b^{4}(M) + \half b^{3}(X) - b^{1}(M) - 1$.
\end{lem}

\begin{proof}
As before we let $E^{2}_{8} = E^{2} \cap H^{2}_{8}(X)$ etc.
As a consequence of corollary \ref{hsuacor} and theorem \ref{absderhamthm}
we find that
$\pi_{H} : \harm^{3}_{cyl} \to H^{3}(M)$ is surjective, and the kernel
is mapped bijectively to $E^{2}_{8}$ by
$\pi_{H,e} : \harm^{3}_{cyl} \to H^{2}(X)$.
Hence $\dim \harm^{3}_{cyl} = b^{3}(M) + \dim E^{2}_{8}$.

The dimension of $E^{2}$ can be computed from the long exact sequence
(\ref{relexacteq}) for relative cohomology together with
(\ref{b3cptdimeq}).
\begin{multline*}
\dim E^{2} = \dim \ker (e : H^{3}_{cpt}(M) \to H^{3}(M)) = 
b^{4}(M) - \dim H^{3}_{0}(M) = b^{4}(M) - b^{3}(M) +\half b^{3}(X)
\end{multline*}
By propositions \ref{aeisoprop} and \ref{cylbochnerprop}
\[ \dim E^{2}_{6} = \dim A^{1} = b^{1}(M),
\;\: \dim E^{2}_{1} = \dim A^{0} = 1 . \]
Hence
\begin{equation*}
\dim E^{2}_{8} =  b^{4}(M) - b^{3}(M) +\half b^{3}(X) - b^{1}(M) - 1 . \qedhere
\end{equation*}
\end{proof}

\section{A topological criterion for $Hol = G_{2}$}
\label{topcondsec}

In this section we prove theorem \ref{maintopcondthm}, which
gives a topological criterion for when the holonomy
group of an EAC \gmfd{} $M^{7}$ is precisely
$G_{2}$ and not a proper subgroup.
As stated in corollary \ref{g2holcor} the holonomy group of a metric
defined by a torsion-free \gstr{} is always a subgroup of~$G_{2}$.
For compact \gmfd s there is a known necessary and sufficient
condition for the holonomy group to be exactly $G_{2}$.

\begin{thm}[{\cite[Proposition $10.2.2$]{joyce00}}] 
\label{cpttopcondthm}
Let $M^{7}$ be a compact \gmfd. Then
$Hol(M) = G_{2}$ if and only if the fundamental group $\pi_{1}(M)$ is finite.
\end{thm}

We summarise the proof, and then generalise the result to the EAC case.
Note that all covering spaces
will be presumed to be equipped with the Riemannian metric pulled back
by the covering map.
In particular all covering maps will be local isometries, and all covering
transformations are isometries.

It is a consequence of the Cheeger-Gromoll line splitting theorem
that any compact Ricci-flat Riemannian manifold $M$ has a finite
cover isometric to a Riemannian product $T^{k} \times N$,
where $T^{k}$ is a flat torus (of dimension $k$ possibly $0$)
and $N$ is compact and simply-connected
(see \cite[Corollary $6.67$]{besse87}).
So for a \gmfd{} $M$ let $\tildm$ be a cover of that form.

If $\pi_{1}(M)$ is infinite then $\tildm = T^{k} \times N$
with $k > 0$, so $Hol(\tildm) \subseteq SU(3)$.
$Hol(\tildm)$ is a finite quotient of $Hol(M)$, so $Hol(M)$
cannot be~$G_{2}$.

If $\pi_{1}(M)$ is finite then $\tildm$ is the universal cover of $M$.
By the classification of Riemannian holonomy groups (`Berger's list',
see e.g. \cite[Theorem $3$]{berger55} or \cite[Theorem $3.4.1$]{joyce00})
the only proper subgroups of $G_{2}$ that can be
the holonomy group of a simply-connected Riemannian
manifold are $1, SU(2)$ and~$SU(3)$ (up to conjugacy).
Thus if $Hol(\tildm)$ is not $G_{2}$ then
it fixes at least one vector in its action on $\bbr^{7}$. By proposition
\ref{simpleholprop} this implies that there exists a parallel $1$-form
$\phi$ on~$\tildm$. Since $\tildm$ is compact and Ricci-flat $\phi$ is
harmonic. But Hodge theory gives
an isomorphism $\harm^{1} \to H^{1}(\tildm)$ between harmonic forms and
de Rham cohomology. $\tildm$ is simply-connected, so $b^{1}(\tildm) = 0$ and
there can be no harmonic $1$-forms on~$\tildm$. Hence
$Hol(M) = Hol(\tildm) = G_{2}$.

\smallskip
To generalise theorem \ref{cpttopcondthm} to EAC
\gmfd s $M$ we use that by proposition \ref{cylbochnerprop} the space
of parallel $1$-forms on $M$ is exactly~$\harm^{1}_{0}$, and
that by corollary \ref{degonehodgecor} the natural map
$\harm^{1}_{0} \to H^{1}(M)$ from bounded harmonic forms to de Rham cohomology
is an isomorphism when $M$ has a single end. Recall from theorem
\ref{salurthm} that a Ricci-flat EAC manifold either has a single end or is
isometric to a product cylinder. We also need the following lemma.

\begin{lem}
\label{topcaseslem}
Let $M$ be a Ricci-flat EAC manifold.
\begin{enumerate}
\item if $M$ has a finite normal cover homeomorphic to a cylinder
then $M$ or a double cover of $M$ is homeomorphic to a cylinder,
\item if $\pi_{1}(M)$ is infinite then $M$ has a finite cover $\tildm$ with
$b^{1}(\tildm) > 0$.
\end{enumerate}
\end{lem}

\begin{proof}
\begin{enumerate}
\item If $\tildm$ is a finite normal cover of $M$ homeomorphic to a cylinder
then it is isometric to a product cylinder $Y \times \bbr$.
$M$ is a quotient of $Y \times \bbr$ by a finite group $A$ of isometries. The
isometries are products of isometries of $Y$ and of $\bbr$ (since they
preserve the set of globally distance minimising geodesics
$\{\{y\} \times \bbr : y \in Y\}$). The elements of $A$ have
finite order, so they must act on the $\bbr$ factor as either the identity or as
reflections. The subgroup $B \subseteq A$ which acts as the identity on
$\bbr$ is either all of~$A$, in which case $M$ is the cylinder
$(Y/B) \times \bbr$, or a normal subgroup of index~$2$, in which case
$(Y/B) \times \bbr$ is a cylindrical double cover of $M$. 

\item $M$ is homotopy equivalent to a compact manifold with boundary,
so $\pi_{1}(M)$ is finitely generated. Since $M$ is complete and has
non-negative Ricci curvature volume comparison arguments show that the volume of
balls in the universal cover of $M$ grows polynomially with the radius,
and this can be used to deduce that $\pi_{1}(M)$ has `polynomial growth'
(see Milnor \cite{milnor68} for details).
A result of Gromov \cite[Main theorem]{gromov81} states that any finitely
generated group with polynomial growth has a nilpotent subgroup of finite index.

So let $G_{0} \subseteq \pi_{1}(M)$ be a nilpotent subgroup of finite index.
$G_{0}$ is soluble, so the derived series $G_{i+1} = [G_{i},G_{i}]$ reaches
$1$. Therefore there is a largest $i$ such that $G_{i} \subseteq \pi_{1}(M)$
has finite index.
Let $\tildm$ be the cover of $M$ corresponding to $G_{i} \subseteq \pi_{1}(M)$.
$G_{i}/G_{i+1}$ is an infinite Abelian group, so has non-zero rank. Hence
\[ b^{1}(\tildm) =
rk \left(\pi_{1}(\tildm)/[\pi_{1}(\tildm),\pi_{1}(\tildm)] \right)
= rk(G_{i}/G_{i+1}) > 0 . \qedhere \]
\end{enumerate}
\end{proof}

\noindent The lemma implies that if $M$ is an EAC \gmfd{} then one of 4
possible cases holds:
\begin{enumerate}

\item \label{cylcase}
$\pi_{1}(M)$ is finite and $M$ is homeomorphic to a cylinder. Then
$M$ is isometric to $Y \times \bbr$ for some compact Calabi-Yau manifold
$Y^{6}$. The arguments in the proof of theorem \ref{cpttopcondthm}
show that the holonomy of $Y$ cannot be a proper subgroup of $SU(3)$.
Thus $Hol(M) = SU(3)$.

\item \label{dcylcase}
$\pi_{1}(M)$ is finite, $M$ has a single end, and has a double cover
homeomorphic to a cylinder. Then the double cover has holonomy $SU(3)$,
so $Hol(M) \not = G_{2}$.

\item \label{infcase}
$\pi_{1}(M)$ is infinite. Then $M$ has a finite cover $\tildm$ with
$b^{1}(M) > 0$. $\tildm$ has a parallel $1$-form by theorem \ref{absderhamthm}
together with proposition \ref{cylbochnerprop}.
Thus $Hol(\tildm) \subseteq SU(3)$, so $Hol(M) \not = G_{2}$.

\item \label{irredcase}
$\pi_{1}(M)$ is finite and neither $M$ nor any double cover of $M$ is
homeomorphic to a cylinder. Then the universal cover $\tildm$ of $M$ is an
EAC \gmfd{} with a single end.
The only proper subgroups of $G_{2}$ that can be the holonomy group of
a complete simply-connected manifold are $1$, $SU(2)$ and $SU(3)$,
so if $Hol(\tildm)$ is not $G_{2}$ then there is a parallel vector field
on $\tildm$. But $b^{1}(\tildm) = 0$, so by corollary \ref{degonehodgecor} and
proposition \ref{cylbochnerprop} there are no parallel $1$-forms on $\tildm$.
Hence $Hol(M) = G_{2}$. 
\end{enumerate}

\smallskip

$Hol(M)$ is exactly $G_{2}$ only in case \ref{irredcase},
so we have proved theorem \ref{maintopcondthm}.
Examples of the cases \ref{dcylcase} and \ref{infcase}
are provided below.
A future paper \cite{jn2} will show that EAC manifolds with holonomy exactly
$G_{2}$ can be constructed by methods similar to those used by Joyce to
construct compact examples in \cite{joyce96-I},
and that both $SU(3)$ and proper subgroups such as $\bbz_{2} \ltimes SU(2)$
can occur as the holonomy of the cross-section of such manifolds.
In particular the converse to the following corollary of theorem
\ref{maintopcondthm} is false.

\begin{cor}
Let $M$ be an EAC \gmfd{} with cross-section $X$, and suppose that $M$
is not finitely covered by a cylinder.
If $Hol(X) = SU(3)$ then $Hol(M) = G_{2}$.
\end{cor}

\begin{proof}
Suppose that $Hol(M)$ is a proper subgroup of $G_{2}$.
Then $\pi_{1}(M)$ is infinite, so $M$ has a finite cover $\tildm$ with
$b^{1}(\tildm) > 0$.
Let $\tilde X$ be the cross-section of $\tildm$. By proposition
\ref{realsubprop} $b^{1}(\tilde X) = 2b^{1}(\tildm) > 0$, so $Hol(\tilde X)$
is a proper subgroup of $SU(3)$. $\tilde X$ is a finite cover of $X$, so
it follows that $Hol(X)$ is a proper subgroup of $SU(3)$.
\end{proof}

\begin{ex}
There exist EAC manifolds $W^{6}$ with holonomy
precisely $SU(3)$ (see Kovalev \cite[Theorem $2.7$]{kovalev03}). Then we can
define a torsion-free \gstr{} on the product $W \times S^{1}$ as in
proposition \ref{cylg2prop}. Of course $Hol(W \times S^{1})$ is not all of
$G_{2}$, but just $SU(3)$. Furthermore \mbox{$b^{1}(W \times S^{1}) > 0$},
so by theorem \ref{maintopcondthm} \emph{no} EAC
\gstr{} on $W \times S^{1}$ can have holonomy exactly~$G_{2}$.
\end{ex}

\begin{ex}
\label{cyinvex}
Let $Y \subset \bbc P^{5}$ be the complex projective variety defined by the
equations $\sum X^{2}_{i} = 0$, $\sum X^{4}_{i} = 0$. $Y$ is a complete
intersection of hypersurfaces, so is a smooth complex 3-fold.
As described in \cite[page 40]{joyce00} the adjunction formula can be used
to show that the first Chern class $c_{1}(Y)$ vanishes.
The Lefschetz hyperplane theorem, stated in the form \cite[Theorem I]{bott59},
can be applied to show that $\pi_{1}(Y) = 1$.

Since the polynomials defining
$Y$ are real the complex conjugation map on $\bbc P^{5}$  restricts to an
involution $a : Y \to Y$. $a$ is anti-holomorphic, and since the
defining polynomials have no roots over $\bbr$ the involution has no
fixed points.

Let $\omega_{FS}$ be the restriction of the Fubini-Study Kähler form to~$Y$.
$a^{*}\omega_{FS} = -\omega_{FS}$. Since $c_{1}(Y) = 0$ Yau's solution to
the Calabi conjecture \cite{yau78} 
implies that there is a unique Kähler form $\omega$ in the cohomology class of
$\omega_{FS}$ such that the corresponding metric is Ricci-flat, making
$Y$ into a Calabi-Yau manifold. The cohomology class of $\omega_{FS}$
is preserved by $-a^{*}$ and $-a^{*}\omega$ is a Kähler form defining a
Ricci-flat metric, so the uniqueness part of the assertion implies that
$-a^{*}\omega = \omega$ (cf. \cite[Proposition $15.2.2$]{joyce00}).

Pick a global holomorphic non-vanishing $3$-form $\phi$ on~$Y$.
$\overline{a^{*}\phi}$ is also holomorphic, so equals $\lambda^{2}\phi$
for some $\lambda \in \bbc$. Replacing $\phi$ with $\lambda\phi$ we can assume
without loss of generality that $\lambda = 1$. Then $\Omega = \re \phi$ is
preserved by $a^{*}$.
We can rescale $\Omega$ to ensure that $(\Omega, \omega)$ is a
\cystr{} in the sense of definition \ref{cy3def}.

Now define a \gstr{} on $Y \times \bbr$ by $\varphi =
\Omega + dt \wedge \omega$. By proposition \ref{cylg2prop}
$\varphi$ is torsion-free.
Let $M$ be the quotient of $Y \times \bbr$ by $a \times (-1)$.
$M$ has a single end and $\pi_{1}(M)$ has order 2.
$(a \times (-1))^{*}\varphi = a^{*}\Omega + (-dt)\wedge a^{*}\omega = \varphi$,
so $\varphi$ induces a well-defined torsion-free \gstr{} on~$M$.
\end{ex}

\emph{Acknowledgements.}
I am grateful to Alexei Kovalev for many helpful discussions,
and to Gabriel Paternain for useful remarks on theorem \ref{maintopcondthm}.

\bibliographystyle{plain}
\bibliography{g2geom}

\end{document}